\documentclass{elsart}
\usepackage{amssymb}
\usepackage[active]{srcltx}

\journal{J. Math. Pures Appl.}

\newcommand{\be}{\begin{equation}}
\newcommand{\ee}{\end{equation}}
\newcommand{\bea}{\begin{eqnarray*}}
\newcommand{\eea}{\end{eqnarray*}}

\newtheorem{theorem}{Theorem}[section]
\newtheorem{corollary}{Corollary}

\newtheorem{proposition}{ Proposition}[section]
\newtheorem{definition}{Definition}[section]
\newtheorem{example}{ Example}[section]

\newtheorem{remark}{ Remark}[section]

\begin{document}
\begin{frontmatter}

\title{Nash-type equilibria on Riemannian manifolds: a variational approach}
\author{Alexandru  Krist\' aly\thanksref{cor}\thanksref{CEU-thanks}}
\address{Department of Economics,  Babe\c s-Bolyai University, 400591 Cluj-Napoca, Romania}

\thanks[cor]{\small Email address: {\rm alexandrukristaly@yahoo.com}}
\thanks[CEU-thanks]{\small Research supported by a grant of the Romanian National Authority for
Scientific Research, CNCS-UEFISCDI, project no.
PN-II-ID-PCE-2011-3-0241, "Symmetries in elliptic problems:
Euclidean and non-Euclidean techniques". The present work was
initiated during the author's visit at Institut des Hautes \'Etudes
Scientifiques (IH\'ES), Bures-sur-Yvette, France.}

 \begin{abstract} {\footnotesize
Motivated by Nash equilibrium problems on 'curved' strategy sets,
the concept of Nash-Stampacchia equilibrium points is introduced via
variational inequalities on Riemannian manifolds. Characterizations,
existence, and stability of Nash-Stampacc\-hia equilibria are
studied when the strategy sets are compact/non\-compact geodesic
convex subsets of Hadamard manifolds, exploiting two well-known
geometrical features of these spaces both involving the metric
projection map. These
 properties actually characterize the {non-positivity} of the
sectional curvature of complete and simply connected Riemannian
spaces, delimiting the Hadamard manifolds as the optimal geometrical
framework of Nash-Stampacc\-hia equilibrium problems. Our analytical
approach exploits various elements from set-valued and variational
analysis, dynamical systems, and non-smooth calculus on Riemannian
manifolds. Examples are presented on the Poincar\'e upper-plane
model and on the open convex cone of symmetric positive definite
matrices endowed with the trace-type Killing form.}\\
\vspace{0.2cm}
\noindent {\bf R\'esum\'e}\\

\noindent    {\footnotesize  Motiv\'e par des probl\`emes
d'\'equilibres de Nash sur des ensembles "courb\'es" de
strat\'egies, la notion d'\'equilibre de Nash-Stampacchia peut
\^etre introduite par des in\'egalité\'es variationnelles sur des
vari\'et\'es Riemanniennes. On \'etudie la caract\'erisa\-tion,
l'existence et la stabilit\'e d'\'equilibres de Nash-Stampacchia
quand les ensembles de strat\'egies sont des sous-ensembles
g\'eod\'esiquement convexes, compacts ou non-compacts, de
vari\'et\'es d'Hadamard, en exploitant deux propri\'et\'es
g\'eom\'etriques bien connues de ces espaces, bas\'ees sur la
projection m\'etrique. En fait ces propri\'et\'es caract\'erisent la
non-positivit\'e de la courbure sectionnelle des espaces de Riemann
complets et simplement connexes, en identifiant les vari\'et\'es
d'Hadamard comme la structure g\'eom\'etrique optimale o\'u poser
les probl\`emes d'\'equilibre de Nash-Stampacchia. Notre approche
analytique utilise
    plusieurs \'el\'ements d'analyse variationnelle et multi-valu\'ee, des syst\`emes
    dynamiques et le calcule non-lisse sur les vari\'et\'es Riemanniennes. Des exemples sont
    pr\'esent\'es dans le cadre du demi-plan de Poincar\'e et dans le c\^one
    ouvert convexe des matrices d\'efinie positives muni d'une forme Killing de
    type trace.}
 \end{abstract}

\tableofcontents


\end{frontmatter}

\section{Introduction}
 After the seminal paper of Nash \cite{N1} there
has been considerable interest in the theory of Nash equilibria due
to its applicability in various real-life phenomena (game theory,
price theory, networks, etc). Appreciating  Nash's contributions, R.
B. Myerson states that "Nash's theory of noncooperative games should
now be recognized as one of the outstanding intellectual advances of
the twentieth century", see also \cite{Mye}. The Nash equilibrium
problem involves $n$ players such that each player know the
equilibrium strategies of the partners, but moving away from his/her
own strategy alone a player has nothing to gain. Formally, if the
sets $K_i$ denote the strategies of the players and
$f_i:K_1\times...\times K_n\to \mathbf R$ are their loss-functions,
$i\in \{1,...,n\},$ the objective is to find an $n$-tuple ${\bf
p}=(p_1,...,p_n)\in {\bf K}=K_1\times...\times K_n$ such that
$$f_i({\bf p})=f_i(p_i,{\bf p}_{-i})\leq f_i(q_i,{\bf p}_{-i})\ {\rm
for\ every}\  q_i\in K_i\ {\rm  and}\ i\in \{1,...,n\},$$  where
$(q_i,{\bf p}_{-i})=(p_1,...,p_{i-1},q_i,p_{i+1},...,p_n)\in{\bf
K}.$ Such point ${\bf p}$ is  called a {\it Nash equilibrium
point for} $({\bf f,K})= (f_1,...,f_n;K_1,...,K_n)$, the set of these points being denoted by ${\mathcal S}_{NE}({\bf f,K}).$  

While most of the known developments in the Nash equilibrium theory
deeply exploit the usual convexity of the sets $K_i$ together with
the vector space structure of their ambient spaces $M_i$ (i.e.,
$K_i\subset M_i$), it is nevertheless true that these results are in
large part {\it geometrical} in nature. The main purpose of this
paper is to enhance those geometrical and analytical structures
which serve as a basis of a systematic study of location of
Nash-type equilibria in a general setting as presently possible. In
the light of these facts our contribution to the Nash equilibrium
theory should be considered intrinsical and analytical rather than
game-theoretical. However, it seems that some ideas of the present
paper can be efficiently applied to evolutionary game dynamics on
curved spaces, see  Bomze \cite{Bomze}, and Hofbauer and Sigmund
\cite{HS}.

Before to start describing our results, we point out an important
(but neglected) topological achievement of Ekeland \cite{Ekeland}
concerning the {\it existence} of Nash-type equilibria for
two-person games on compact manifolds based on transversality and
fixed point arguments. Without the sake of completeness, Ekeland's
result says that if $f_1,f_2:M_1\times M_2\to \mathbf R$ are
continuous functions having also some differentiability properties
where $M_1$ and $M_2$ are compact manifolds whose Euler-Poincar\'e
characteristics are non-zero (orientable case) or odd
(non-orientable case), then there exists at least a Nash-type
equilibria for $(f_1,f_2;M_1,M_2)$ formulated via first order
conditions involving the terms $\frac{\partial f_i}{\partial x_i}$,
$i=1,2.$

In the present paper we assume {\it a priori} that the strategy sets
$K_i$ are {\it geodesic convex} subsets of certain
finite-dimensional Riemannian manifolds $(M_i,g_i)$, i.e., for any
two points of $K_i$ there exists a unique geodesic in $(M_i,g_i)$
connecting them which belongs entirely to $K_i$. This approach can
be widely applied when the strategy sets are 'curved'. Note that the
choice of such Riemannian structures does not influence the Nash
equilibrium points for $({\bf f,K})$. As far as we know, the first
step into this direction was made recently in
\cite{Kristaly-PAMS-2009} via a McClendon-type minimax inequality
for acyclic ANRs, guaranteeing the existence of at least one Nash
equilibrium point for $({\bf f,K})$ whenever $K_i\subset M_i$ are
compact and geodesic convex sets of certain finite-dimensional
Riemannian manifolds $(M_i,g_i)$ while the functions $f_i$ have
certain regularity properties, $i\in \{1,...,n\}$. By using
Clarke-calculus on manifolds, in \cite{Kristaly-PAMS-2009} we
introduced and studied for a wide class of {\it non}-smooth
functions the set of {\it Nash-Clarke equilibrium points for} ${ \bf
(f,K)}$, denoted in the sequel as
${\mathcal S}_{NC}({\bf f,K})$; 
 see Section \ref{exist-Nash}. Note that ${\mathcal
S}_{NC}({\bf f,K})$ is larger than ${\mathcal S}_{NE}({\bf f,K})$;
thus, a promising way to localize the elements of ${\mathcal
S}_{NE}({\bf f,K})$ is to determine the set ${\mathcal S}_{NC}({\bf
f,K})$ and to choose among these points the appropriate ones. In
spite of the naturalness of this approach, we already pointed out
its limited applicability due to the involved definition of
${\mathcal S}_{NC}({\bf f,K})$, conjecturing a more appropriate
concept in order to locate the elements of ${\mathcal S}_{NE}({\bf
f,K})$.

Motivated by the latter problem, we observe that the Fr\'echet and
limiting sub\-differential calculus of lower semicontinuous
functions on Riemannian manifolds developed by Azagra, Ferrera and
L\'opez-Mesas \cite{Azagra-JFA2},  and Ledyaev and Zhu
\cite{LZ-TAMS} provides a satisfactory approach. The idea is to
consider the following {\it system of variational inequalities}:
find ${\bf p}\in {\bf K}$ and $\xi_C^i\in
\partial_C^i f_i({\bf p})$ such that
$$\langle \xi_C^i, \exp_{p_i}^{-1}(q_i)\rangle_{g_i} \geq 0  \
{\rm for\ all}\ q_i\in K_i,\ i\in \{1,... ,n\},$$ where
$\partial_C^i f_i({\bf p})$ denotes the Clarke subdifferential of
the locally Lipschitz function $f_i(\cdot,{\bf p}_{-i})$ at the
point $p_i\in K_i;$ for details, see Section \ref{exist-Nash}.
 The solutions of this system form the set of {\it Nash-Stampacchia\footnote{Terminology introduced in honor of G. Stampacchia for his
 deep contributions to the theory of variational inequalities.} equilibrium points for} ${ \bf
(f,K)}$, denoted by ${\mathcal S}_{NS}({\bf f,K}),$ which is the
main concept of the present paper.

Our first result shows that $${\mathcal S}_{NE}({\bf f,K})\subset
{\mathcal S}_{NS}({\bf f,K})= {\mathcal S}_{NC}({\bf f,K})$$ for the
same class of non-smooth functions ${\bf f}=(f_1,...,f_n)$ as in
\cite{Kristaly-PAMS-2009} (see Theorem \ref{equi-kritikus} (i)).
Although ${\mathcal S}_{NS}({\bf f,K})$ and ${\mathcal S}_{NC}({\bf
f,K})$ coincide, the set of Nash-Stampacchia equilibrium points is
more flexible and applicable. Indeed, via this new notion we can
handle both compact and non-compact strategy sets, and we identify
the geometric framework where this argument works, i.e., the class
of Hadamard manifolds. These kinds of results could not have been
achieved via the notion of  Nash-Clarke equilibrium points as we
will describe later.

To establish the above inclusions we introduce a  notion of
subdifferential for non-smooth functions by means of the cut locus
which could be of some interest in its own right as well. Then,
explicit characterizations of the Fr\'echet and limiting normal
cones of geodesic convex sets in arbitrarily Riemannian manifolds
are given by exploiting some fundamental results from
\cite{Azagra-JFA2} and \cite{LZ-TAMS}. If ${\bf f}=(f_1,...,f_n)$
verifies a suitable 'diagonal' convexity assumption then we have
equalities in the above relation (see Theorem \ref{equi-kritikus}
(ii)).

Having these inclusions in mind, the main purpose of the present
paper is to establish existence, location and stability of
Nash-Stampacchia equilibrium points for ${ \bf (f,K)}$ in different
settings. While a Nash equilibrium point is obtained precisely as
the fixed point of a suitable function (see for instance Nash's
original proof  via Kakutani fixed-point theorem), Nash-Stampacchia
equilibrium points are expected to be characterized in a similar way
as fixed points of a special map defined on the product Riemannian
manifold ${\bf M}=M_1\times...\times M_n$ endowed with its natural
Riemannian metric ${\bf g}$ inherited from the metrics $g_i$, $i\in
\{1,...,n\}$. In order to achieve this aim, certain curvature and
topological restrictions are needed on the manifolds $(M_i,g_i)$. By
assuming that the ambient Riemannian manifolds $(M_i,g_i)$ for the
geodesic convex strategy sets $K_i$ are {\it Hadamard manifolds},
our key observation (see Theorem \ref{theorem-equivalence}) is that
${\bf p\in K}$ is a Nash-Stampacchia equilibrium point for $({\bf
f,K})$ {if and only if} ${\bf p}$ is a fixed point of the set-valued
map $A_\alpha^{\bf f}:{\bf K}\to 2^{\bf K}$ defined by
$$A_\alpha^{\bf f}({\bf p})=P_{\bf K}(\exp_{\bf
p}(-\alpha\partial_C^\Delta {\bf f}({\bf p}))).$$ Here, $P_{\bf K}$
is the metric projection operator associated to the geodesic convex
set ${\bf K\subset M}$, $\alpha>0$ is a fixed number, and
$\partial_C^\Delta {\bf f}({\bf p})$ denotes the diagonal Clarke
subdifferential at point ${\bf p}$ of ${\bf f}=(f_1,...,f_n)$; see
Section \ref{exist-Nash}.

Within this geometrical framework, two cases are discussed. On the
one hand, when ${\bf K\subset M}$ is {\it compact}, one can prove
via the Begle's fixed point theorem for set-valued maps the
existence of at least one Nash-Stampacchia equilibrium point for
$({\bf f,K})$ (see Theorem \ref{begle}). On the other hand, when
${\bf K\subset M}$ is {\it not} necessarily {\it compact}, we
provide two types of results. First, based on a suitable coercivity
assumption on $\partial_C^\Delta {\bf f}$, combined with the result
for compact sets, we are able to guarantee the existence of at least
one Nash-Stampacchia equilibrium point for $({\bf f,K})$ (see
Theorem \ref{begle-noncompact}). Second, by requiring more
regularity on ${\bf f}$ in order to avoid technicalities, we
consider two dynamical systems; a discrete one
$$\leqno{(DDS)_\alpha}\ \ \ \ \ {\bf p}_{k+1}= A_\alpha^{\bf
f}(P_{\bf K}({\bf p}_k)),\ \ \ {\bf p}_0\in {\bf M};$$ and a
continuous one
$$\leqno{(CDS)_\alpha} \ \ \ \ \ \left\{
\begin{array}{lll}
\dot\eta(t)= \exp_{\eta(t)}^{-1}(A_\alpha^{\bf f}(P_{\bf
K}(\eta(t))))
\\ \eta(0)={\bf
p}_0\in {\bf M}.
\end{array}
\right. $$ By assuming a Lipschitz-type condition on
$\partial_C^\Delta {\bf f}$, one can  prove that the set of
Nash-Stampacchia equilibrium points for $({\bf f,K})$ is a {\it
singleton} and the orbits of both dynamical systems exponentially
converge to this unique point (see Theorem
\ref{dyn-system-theorem-fo}). Here, we exploit some arguments from
the theory of differential equations on manifolds as well as careful
comparison results of Rauch-type. It is clear by construction that
the orbit of $(DDS)_\alpha$ is viable relative to the set ${\bf K}$,
i.e., ${\bf p}_k\in {\bf K}$ for every $k\geq 1$. By using a recent
result of Ledyaev and Zhu \cite{LZ-TAMS}, one can also prove an
invariance property of the set ${\bf K}$ with respect to the orbit
of $(CDS)_\alpha$.
Note that the aforementioned results concerning the 'projected'
dynamical system $(CDS)_\alpha$ are new even in the Euclidean
setting; see Cavazzuti, Pappalardo and Passacantando
\cite{CPP-JOTA}, Xia \cite{Xia-JOTA}, and Xia and  Wang
\cite{XW-JOTA}.

Since the manifolds $(M_i,g_i)$ are assumed to be of Hadamard type
(see Theorems
\ref{theorem-equivalence}-\ref{dyn-system-theorem-fo}), so is the
product manifold ${\bf (M,g)}$. Our analytical arguments
deeply exploit two geometrical features of the product {\it Hadamard
manifold} ${\bf (M,g)}$ concerning the metric projection operator
for closed, geodesic convex sets:
\begin{itemize}
  \item [(A)] {\it Validity of the obtuse-angle property}, see Proposition
\ref{egyik-irany} (i). This fact is exploited in the
characterization of Nash-Stampacchia equilibrium points for $({\bf
f,K})$ via the fixed points of the map $A_\alpha^{\bf f}$, see
Theorem \ref{theorem-equivalence}.

  \item [(B)] {\it  Non-expansiveness of the projection operator}, see
Proposition \ref{egyik-irany} (ii). This property is applied several
times in the proof of Theorems \ref{begle} and
\ref{dyn-system-theorem-fo}.
\end{itemize}
\noindent It is natural to ask to what extent the Riemannian
structures of $(M_i,g_i)$ are determined when the properties (A) and
(B) simultaneously hold on the product manifold ${\bf (M,g)}$. A
constructive proof combined with the formula of sectional curvature
via the Levi-Civita parallelogramoid and a result of Chen
\cite{CH-Chen-TAMS}  shows that if $(M_i,g_i)$ are complete, simply
connected Riemannian manifolds then (A) and (B) are simultaneously
verified on ${\bf (M,g)}$ if and only if $(M_i,g_i)$ are Hadamard
manifolds (see Theorem \ref{fotetel-1}). Consequently, we may assert
that {Hadamard manifolds} are the optimal geometrical framework to
elaborate a fruitful theory of Nash-Stampacchia equilibrium problems
on Riemannian manifolds. Furthermore, we notice that properties (A)
and (B) are also the milestones of the theory of monotone vector
fields, proximal point algorithms and variational inequalities
developed on Hadamard manifolds by Li, L\'opez and Mart\'\i
n-M\'arquez \cite{LLMM-1}, \cite{LLMM-2}, and N\'emeth
\cite{Nemeth}. As a byproduct of Theorem \ref{fotetel-1} we state
that Hadamard manifolds are the appropriate geometrical frameworks
among Riemannian manifolds for the aforementioned theories.

The paper is divided as follows. In \S \ref{preli-sect} we
recall/prove those notions and results which will be used throughout
the paper: basic elements from Riemannian geometry, the
parallelogramoid of Levi-Civita; properties of the metric
projection; non-smooth calculus, dynamical systems and viability
results on Riemannian manifolds. In \S \ref{exist-Nash} we compare
the three Nash-type equilibria; simultaneously, we also recall some
results from \cite{Kristaly-PAMS-2009}.  In \S \ref{section-main},
we prove the main results of this paper both for compact and
non-compact strategy sets which are 'embedded' into certain Hadamard
manifolds.  In \S \ref{section-Nash-curvature} we characterize the
geometric properties (A) and (B) on ${\bf (M,g)}$ by the Hadamard
structures of the complete and simply connected Riemannian manifolds
$(M_i,g_i)$, $i\in \{1,...,n\}$. Finally, in \S \ref{utol-sect} we
present some relevant examples on the Poincar\'e upper-plane model
as well as on the Hadamard manifold of symmetric positive definite
matrices endowed with the Killing form of trace-type. Our examples
are motivated by some applications from Bento, Ferreira and Oliveira
\cite{BFO-1}, \cite{BFO-2}, Colao, L\'opez, Marino and
Mart\'in-M\'arquez \cite{CLMMM}, and Li and Yao \cite{Li-Yao}.

\section{Preliminaries: metric projections, non-smooth calculus and dynamical systems on Riemannian manifolds}\label{preli-sect}

2.1. {\bf Elements from Riemannian geometry.}
We first recall those elements from Riemannian geometry which will
be used throughout the paper. We mainly follow Cartan \cite{Cartan}
and do Carmo \cite{doCarmo}.

In this subsection, $(M,g)$ is a connected $m$-dimensional
Riemannian manifold, $m\geq 2$. Let $TM=\cup_{p\in M}(p,T_pM)$ and
$T^*M=\cup_{p\in M}(p,T_p^*M)$ be the tangent and cotangent bundles
to $M.$ For every $p\in M$, the Riemannian metric induces a natural
Riesz-type isomorphism between the tangent space $T_pM$ and its dual
$T_p^*M$; in particular, if $\xi\in T_p^*M$ then there exists a
unique $W_\xi\in T_pM$ such that
\begin{equation}\label{dual-nem-dual}
\langle\xi,V\rangle_{g,p}=g_p(W_\xi,V)\ \mbox{for all}\ V\in T_pM.
\end{equation}
Instead of $g_p(W_\xi,V)$ and $\langle\xi,V\rangle_{g,p}$ we shall
write simply $g(W_\xi,V)$ and $\langle\xi,V\rangle_g$ when no
confusion arises. Due to (\ref{dual-nem-dual}), the elements $\xi$
and $W_\xi$ are identified. With the above notations, the norms on
$T_pM$ and $T_p^*M$ are defined by
$$\|\xi\|_g=\|W_\xi\|_g=\sqrt{g(W_\xi,W_\xi)}.$$ The
generalized Cauchy-Schwartz inequality is also valid, i.e.,  for
every $V\in T_pM$ and $\xi\in T_p^*M$,
\begin{equation}\label{cauchy}
  |\langle \xi, V\rangle_g|\leq \|\xi\|_g\|V\|_g.
\end{equation}
Let $\xi_k\in T_{p_k}^*M$, $k\in \mathbf N$, and $\xi\in T_p^*M$.
The sequence $\{\xi_k\}$ converges to $\xi$, denoted by $\lim_k
\xi_k=\xi$, when $p_k\to p$ and $\langle \xi_k,
W(p_k)\rangle_g\to\langle \xi, W(p)\rangle_g$ as $k\to \infty$, for
every $C^\infty$ vector field $W$ on $M$.

Let $h:M\to \mathbf R$ be a $C^1$ functional at $p\in M$; the
differential of $h$ at $p$, denoted by $dh(p)$, belongs to $T_p^*M$
and is defined by $$\langle dh(p),V\rangle_{g}=g({\rm grad} h(p),V)\
\mbox{for all}\ V\in T_pM.$$ If $(x^1,...,x^m)$ is the local
coordinate system on a coordinate neighborhood $(U_p,\psi)$ of $p\in
M$,  and the  local components of $dh$ are denoted
$h_i=\frac{\partial h}{\partial x_i}$, then the local components of
grad$h$ are $h^i=g^{ij}h_j$. Here, $g^{ij}$ are the local components
of $g^{-1}$.

Let $\gamma:[0,r]\to M$ be a $C^1$ path, $r>0$. The length of
$\gamma$ is defined by $$L_g(\gamma)=\int_0^r\|\dot
\gamma(t)\|_gdt.$$ For any two points $p,q\in M$, let
$$d_g(p,q)=\inf\{L_g(\gamma):\gamma\ \mbox{is a}\ C^1\ \mbox{path
joining}\ p\ \mbox{and}\ q\ \mbox{in}\ M\}.$$ The function
$d_g:M\times M\to \mathbf R$ is a metric which generates the same
topology on $M$ as the underlying manifold topology. For every $p\in
M$ and $r>0$, we define the open ball of center $p\in M$ and radius
$r>0$ by $$B_g(p,r)=\{q\in M:d_g(p,q)<r\}.$$

Let us denote by $\nabla$ the unique natural covariant derivative on
$(M,g)$, also called the Levi-Civita connection. A vector field $W$
along a $C^1$ path $\gamma$ is called parallel when $\nabla_{\dot
\gamma}W=0$. A $C^\infty$ parameterized path $\gamma$ is a geodesic
in $(M,g)$ if its tangent $\dot \gamma$ is parallel along itself,
i.e., $\nabla_{\dot\gamma}\dot\gamma=0$. The geodesic segment
$\gamma:[a,b]\to M$ is called minimizing if
$L_g(\gamma)=d_g(\gamma(a),\gamma(b)).$

  Standard ODE theory implies
that for every $V\in T_pM$, $p\in M$, there exists an open interval
$I_V\ni 0$ and a unique geodesic $\gamma_V:I_V\to M$ with
$\gamma_V(0)=p$ and $\dot \gamma_V(0)=V.$ Due to the 'homogeneity'
property of the geodesics (see \cite[p. 64]{doCarmo}), we may define
the exponential map $\exp_p:T_pM\to M$ as $\exp_p(V)=\gamma_V(1).$
Moreover,
\begin{equation}\label{exp-identity}
  d\exp_p(0)={\rm id}_{T_pM}.
\end{equation}
 Note that there exists an open (starlike) neighborhood $\mathcal U$ of the zero vectors in $TM$ and an open neighborhood $\mathcal V$
 of the diagonal $M\times M$ such that the exponential map $V\mapsto \exp_{\pi(V)}(V)$ is smooth and the map
 $\pi\times \exp:\mathcal U\to \mathcal V$ is a diffeomorphism, where $\pi$ is the canonical projection of $TM$ onto $M.$
 Moreover, for every $p\in M$ there exists a number $r_p>0$ and a neighborhood $\tilde U_p$ such
 that for every $q\in \tilde U_p$, the map $\exp_q$ is a $C^\infty$ diffeomorphism on $B(0,r_p)\subset T_qM$ and $\tilde U_p\subset
 \exp_q(B(0,r_p))$; the set $\tilde U_p$ is called a {\it totally normal
 neighborhood} of $p\in M$. In particular, it follows that every two points $q_1,q_2\in
 \tilde U_p$ can be joined by a minimizing geodesic of
 length less than $r_p$. Moreover, for every
$q_1,q_2\in \tilde U_p$ we have
\begin{equation}\label{exp-dist}
\|\exp_{q_1}^{-1}(q_2)\|_g=d_g(q_1,q_2).
\end{equation}
The {\it tangent cut locus} of $p\in M$ in $T_pM$ is the set of all
vectors $v\in T_pM$ such that $\gamma(t) = \exp_p(tv)$ is a
minimizing geodesic for $t \in [0,1]$ but fails to be minimizing for
$t \in [0,1 + \varepsilon)$ for each $\varepsilon > 0.$ The {\it cut
locus of} $p\in M$, denoted by $C_p$,  is the image of the tangent
cut locus of $p$ via $\exp_p$. Note that any totally normal
 neighborhood of $p\in M$ is contained into $M\setminus C_p.$

We conclude this subsection by recalling a less used form of the
{\it sectional curvature} by the so-called Levi-Civita
parallelogramoid. Let $p\in M$ and $V_0,W_0\in T_pM$ two vectors
with $g(V_0,W_0)=0$. Let $\sigma:[-\delta,2\delta]\to M$ be the
geodesic segment $\sigma(t)=\exp_p(tV_0)$ and $W$ be the unique
parallel vector field along $\sigma$ with the initial data
$W(0)=W_0$, the number $\delta>0$ being small enough. For any $t\in
[0,\delta]$, let $\gamma_t:[0,\delta]\to M$ be the geodesic
$\gamma_t(u)=\exp_{\sigma(t)}(uW(t)).$ The sectional curvature of
the subspace $S=$span$\{W_0,V_0\}\subset T_pM$ at the point $p\in M$
is given by $$
  K_p(S)=\lim_{u,t\to 0}\frac{d_g^2(p,\sigma(t))-d_g^2(
\gamma_0(u), \gamma_t(u))}{d_g( p, \gamma_0(u))\cdot d_g( p,
\sigma(t))}, $$ see Cartan \cite[p. 244-245]{Cartan}. The
infinitesimal geometrical object determined by the four points $p$,
$\sigma(t)$, $\gamma_0(u),$ $\gamma_t(u)$ (with $t,u$ small enough)
is called the {parallelogramoid of Levi-Civita.}\\

\vspace{0.1cm}

\noindent 2.2. {\bf Metric projections.} 
Let $(M,g)$ be an $m$-dimensional Riemannian manifold $(m\geq 2)$,
$K\subset M$ be a non-empty set. Let $$P_K(q)=\{p\in
K:d_g(q,p)=\inf_{z\in K}d_g(q,z)\}$$ be the set of {\it metric
projections} of the point $q\in M$ to the set $K$. Due to the
theorem of Hopf-Rinow, if $(M,g)$ is complete, then any closed set
$K\subset M$ is {\it proximinal}, i.e., $P_K(q)\neq \emptyset$ for
all $q\in M$. In general, $P_K$ is a set-valued map. When $P_K(q)$
is a singleton for every $q\in M,$ we say that $K$ is a {\it
Chebyshev set}. The map $P_K$ is {\it non-expansive} if
$$d_g(p_1,p_2)\leq d_g(q_1,q_2)\ \ {\rm for\ all}\
q_1,q_2\in M\ {\rm and}\ p_1\in P_K(q_1), p_2\in P_K(q_2).$$ In
particular, $K$ is a Chebyshev set whenever the map $P_K$ is
non-expansive.

The set $K\subset M $ is {\it geodesic convex} if every two points
$q_1,q_2\in K$ can be joined by a unique minimizing geodesic whose
image belongs
 to $K.$ Note that (\ref{exp-dist}) is also valid for every $q_1,q_2\in K$
in a geodesic convex set $K$ since $\exp_{q_i}^{-1}$ is well-defined
on $K$, $i\in \{1,2\}$. The function $f:K\to \mathbf R$ is {\it
convex}, if $f\circ \gamma:[0,1]\to \mathbf R$ is convex in the
usual sense for every geodesic $\gamma:[0,1]\to K$ provided that
$K\subset M$ is a geodesic convex set.

 A non-empty closed set $K\subset M$ verifies the {\it obtuse-angle property} if for
fixed $q\in M$ and $p\in K$ the following two statements are
equivalent: \vspace{0.2cm}
\begin{itemize}
\item[$(OA_1)$]  $p\in P_K(q)$;
\item[$(OA_2)$]  If $\gamma:[0,1]\to M$ is the unique
minimal geodesic from $\gamma(0)=p\in K$ to $\gamma(1)=q$, then for
every geodesic $\sigma:[0,\delta]\to K$ $(\delta\geq 0)$ emanating
from the point $p$, we have $g(\dot\gamma(0),\dot\sigma(0))\leq 0.$
\end{itemize}

\begin{remark}\label{remark-obtuse-kicsi}\rm
(a) In the Euclidean case $(\mathbf R^m,\langle \cdot,\cdot
\rangle_{\mathbf R^m}),$ (here, $\langle\cdot,\cdot\rangle_{\mathbf
R^m}$ is the standard inner product in $\mathbf R^m$), every
non-empty closed convex set $K\subset \mathbf R^m$ verifies the
obtuse-angle property, see Moskovitz-Dines \cite{MD}, which reduces
to the well-known geometric form: $$p\in P_K(q) \Leftrightarrow
\langle q-p, z-p\rangle_{\mathbf R^m} \leq 0\ \ \mbox{for all}\ z\in
K.$$

(b) The first variational formula of Riemannian geometry shows that
$(OA_1)$ implies $(OA_2)$ for every closed set $K\subset M$ in a
complete Riemannian manifold $(M,g)$. However, the converse does not
hold in general; for a detailed discussion,  see Krist\'aly, R\u
adulescu and Varga \cite{KRV}.

\end{remark}


A Riemannian manifold $(M,g)$ is a {\it Hadamard manifold} if it is
complete, simply connected and its sectional curvature is
non-positive. It is well-known that on a Hadamard manifold $(M,g)$
 every geodesic convex set is a Chebyshev set, see Jost \cite{Jost}. Moreover, we
have

\begin{proposition}\label{egyik-irany} {\it Let $(M,g)$ be a
finite-dimensional Hadamard manifold, $K\subset M$ be a closed set.
The following statements hold true:
 \begin{enumerate}
 \item[{\rm (i)}] $(${\rm Walter \cite{Walter}}$)$  If  $K\subset M$ is
geodesic convex, it verifies the obtuse-angle property;
 \item[{\rm (ii)}] $(${\rm Grognet \cite{Grognet}}$)$ $P_K$ is  non-expansive if and only if $K\subset M$ is geodesic
convex.
 \end{enumerate}}
\end{proposition}
\noindent Finally, we recall that on a Hadamard manifold $(M,g)$, if
$h(p)=d_g^2(p,p_0)$,  $p_0\in M$ is fixed, then
\begin{equation}\label{deriv-tavolsag-hadamard}
    {\rm grad}h(p) =-2\exp_p^{-1}(p_0).
\end{equation}

\vspace{0.3cm}

\noindent 2.3. {\bf Non-smooth calculus on manifolds.} We first
recall some basic notions and results from the subdifferential
calculus on Riemannian manifolds, developed by Azagra, Ferrera and
L\'opez-Mesas \cite{Azagra-JFA2}, Ledyaev and Zhu \cite{LZ-TAMS}.
Simultaneously, we introduce two subdifferential notions based on
the cut locus, and we establish an analytical characterization of
the limiting/Fr\'echet normal cone on Riemannian manifolds (see
Corollary \ref{th-normal-cones-tetel}) which plays a crucial role in
the study of Nash-Stampacchia equilibrium points.

Let $(M,g)$ be an $m$-dimensional Riemannian manifold and let
$f:M\to \mathbf R\cup \{+\infty\}$ be a lower semicontinuous
function with dom$(f)\neq \emptyset.$ The {\it
Fr\'echet-sub\-dif\-fe\-ren\-tial} of $f$ at $p\in {\rm dom}(f)$ is
the set
$$\partial_F f(p)=\{dh(p): h\in C^1(M)\ {\rm and}\ f-h\ {\rm
attains\ a\ local\ minimum\ at}\ p\}.$$

\begin{proposition}{\rm \cite[Theorem 4.3]{Azagra-JFA2}}\label{subdiff-char}
{\it Let $(M,g)$ be an $m$-dimensional Riemannian manifold and let
$f:M\to \mathbf R\cup \{+\infty\}$ be a lower semicontinuous
function,  $p\in${\rm dom}$(f)\neq \emptyset$ and $\xi\in T_p^*M.$
The following statements are equivalent:
\begin{itemize}
\item[{\rm (i)}] $\xi\in \partial_F f(p)$;
\item[{\rm (ii)}] For every chart $\psi:U_p\subset M\to \mathbf R^m$ with $p\in
U_p$, if $\zeta=\xi\circ d\psi^{-1}(\psi(p))$, we have that
 $$\liminf_{v\to 0}\frac{(f\circ \psi^{-1})(\psi(p)+v)-f(p)-\langle \zeta,v\rangle_g}{\|v\|}\geq
 0;$$
\item[{\rm (iii)}] There exists a chart $\psi:U_p\subset M\to \mathbf R^m$ with $p\in
U_p$, if $\zeta=\xi\circ d\psi^{-1}(\psi(p))$, then
$$\liminf_{v\to 0}\frac{(f\circ \psi^{-1})(\psi(p)+v)-f(p)-\langle
\zeta,v\rangle_g}{\|v\|}\geq
 0.$$
\end{itemize}
In addition, if $f$ is locally bounded from below, i.e., for every
$q\in M$ there exists a neighborhood $U_q$ of $q$ such that $f$ is
bounded from below on $U_q$, the above conditions are also
equivalent to
\begin{itemize}
\item[{\rm (iv)}] There exists a function $h\in C^1(M)$ such that  $f-h$
attains a global minimum at $p$ and $\xi=dh(p).$
\end{itemize}}
\end{proposition}

\noindent The {\it limiting subdifferential} and {\it singular
subdifferential of $f$ at} $p\in M$ are the sets $$\partial_L
f(p)=\{\lim_k \xi_k: \xi_k\in
\partial_F f(p_k),\ (p_k,f(p_k))\to (p,f(p))\}$$
and $$\partial_\infty f(p)=\{\lim_k t_k\xi_k: \xi_k\in
\partial_F f(p_k),\ (p_k,f(p_k))\to (p,f(p)), t_k\to 0^+\}.$$

\begin{proposition}\label{prop-lok-minim}{\rm \cite{LZ-TAMS}}
{\it Let $(M,g)$ be a finite-dimensional Riemannian manifold and let
$f:M\to \mathbf R\cup \{+\infty\}$ be a lower semicontinuous
function. Then, we have
\begin{itemize}
  \item[{\rm (i)}] $\partial_F f(p)\subset \partial_L f(p),\ p\in{\rm
  dom}(f);$
  \item[{\rm (ii)}] $0\in \partial_\infty f(p),$ $p\in M;$
  \item[{\rm (iii)}] If $p\in${\rm dom}$(f)$ is a local minimum of $f$, then
  $0\in \partial_F f(p)\subset \partial_L f(p).$
\end{itemize}}
\end{proposition}

\begin{proposition}\label{mean-value}{\rm \cite[Theorem 4.8 (Mean Value inequality)]{LZ-TAMS}} 
Let $f:M\to \mathbf R$ be a continuous function bounded from below,
let $V$ be a $C^\infty$ vector field on $M$ and let $c:[0,1]\to M$
be a curve such that $\dot c(t)=V(c(t))$, $t\in [0,1]$. Then for any
$r<f(c(1))-f(c(0))$, any $\varepsilon>0$ and any open neighborhood
$U$ of $c([0,1])$, there exists $m\in U$, $\xi\in \partial_F f(m)$
such that $r<\langle \xi, V(m)\rangle_g$.
\end{proposition}

\begin{proposition}\label{sum-rule} {\rm \cite[Theorem 4.13 (Sum rule)]{LZ-TAMS}}
{\it Let $(M,g)$ be an $m$-dimensional Riemannian manifold and let
$f_1,...,f_H:M\to \mathbf R\cup \{+\infty\}$ be lower semicontinuous
functions. Then, for every $p\in M$ we have either
$\partial_L(\sum_{l=1}^H  f_l)(p)\subset \sum_{l=1}^H
\partial_L f_l(p)$,
or there exist $\xi_l^\infty\in \partial_\infty f_l(p),$
$l=1,...,H,$ not all zero such that $\sum_{l=1}^H \xi_l^\infty =0.$}
\end{proposition}

The {\it cut-locus subdifferential} of $f$ at $p\in {\rm dom}(f)$ is
defined as
$$\partial_{cl}f(p)=\{\xi\in T_p^*M:f(q)-f(p)\geq \langle \xi,\exp_p^{-1}(q)\rangle_g\ {\rm for\ all}\ q\in M\setminus C_p\},$$
where $C_p$ is the cut locus of the point $p\in M$. Note that
$M\setminus C_p$ is the maximal open set in $M$ such that every
element from it can be joined to $p$ by exactly one minimizing
geodesic, see Klingenberg \cite[Theorem 2.1.14]{Kling}. Therefore,
the cut-locus subdifferential is well-defined, i.e.,
$\exp_p^{-1}(q)$ makes sense and is unique for every $q\in
M\setminus C_p$. We first prove

\begin{theorem}\label{subdifferentials-theorem-uj}
Let $(M,g)$ be a Riemannian manifold and $f : M \to \mathbf R\cup
\{+\infty\}$ be a proper, lower semicontinuous function. Then, for
every $p\in {\rm dom}(f)$  we have
$$\partial_{cl}f(p)\subset \partial_Ff(p) \subset
\partial_Lf(p).$$ Moreover, if $f$ is convex, the above inclusions
become equalities.
\end{theorem}

\noindent {\it Proof.} The  last inclusion is standard, see
Proposition \ref{prop-lok-minim}(i).
  Now, let $\xi\in
\partial_{cl}f(p)$, i.e., $f(q)-f(p)\geq \langle \xi,\exp_p^{-1}(q)\rangle_g\ {\rm for\ all}\ q\in M\setminus
C_p.$ In particular, the latter inequality is valid for every $q\in
B_g(p,r)$ for $r>0$ small enough, since $ B_g(p,r) \subset
M\setminus C_p$  (for instance, when $B_g(p,r)\subset M$ is a
totally normal ball around $p$). Now, by choosing
$\psi=\exp_{p}^{-1}:B_g(p,r)\to T_{p}M$ in Proposition
\ref{subdiff-char}(ii), one has that $f(\exp_pv)-f(p)\geq \langle
\xi,v\rangle_g\ {\rm for\ all}\ v\in T_pM,$ $\|v\|<r,$ which implies
 $\xi\in \partial_Ff(p)$.

Now, we assume in addition that $f$ is convex, and let $\xi\in
\partial_Lf(p).$ We are going to prove that $\xi\in
\partial_{cl}f(p).$ Since $\xi\in
\partial_Lf(p)$, we have that $\xi=\lim_k \xi_k$ where $\xi_k\in
\partial_F f(p_k),\ (p_k,f(p_k))\to (p,f(p))$. By Proposition
\ref{subdiff-char}(ii), for $\psi_k=\exp_{p_k}^{-1}:\tilde
U_{p_k}\to T_{p_k}M$ where $\tilde U_{p_k}\subset M$ is a totally
normal ball centered at $p$, one has that
\begin{equation}\label{equation-cut-locus}
\liminf_{v\to 0}\frac{f(\exp_{p_k}v)-f(p_k)-\langle
\xi_k,v\rangle_g}{\|v\|}\geq
 0.
\end{equation}
Now, fix $q\in M\setminus C_p$. The latter fact is equivalent to
$p\in M\setminus C_q$, see Klingenberg \cite[Lemma 2.1.11]{Kling}.
Since $M\setminus C_q$ is open and $p_k\to p$, we may assume that
$p_k\in M\setminus C_q$, i.e., $q$  and every point $p_k$ is joined
by a unique minimizing geodesic. Therefore, $V_k=\exp_{p_k}^{-1}(q)$
is well-defined. Now, let $\gamma_k(t)=\exp_{p_k}(t V_k)$ be the
geodesic which joins $p_k$ and $q$. Then (\ref{equation-cut-locus})
implies that
\begin{equation}\label{kell-hatrabbb}
   \liminf_{t\to 0^+}\frac{f(\gamma_k(t))-f(p_k)-\langle \xi_k,t
V_k\rangle_g}{\|tV_k\|}\geq
 0.
\end{equation}
Since $f$ is convex, one has that $f(\gamma_k(t))\leq t
f(\gamma_k(1))+(1-t)f(\gamma_k(0)),\ t\in [0,1],$ thus, the latter
relations imply that
$$\frac{f(q)-f(p_k)-\langle \xi_k,
\exp_{p_k}^{-1}(q)\rangle_g}{d_g(p_k,q)}\geq
 0.$$
Since $f(p_k)\to f(p)$ and $\xi=\lim_k \xi_k$, it yields precisely
that $$f(q)-f(p)-\langle \xi, \exp_{p}^{-1}(q)\rangle_g\geq
 0,$$
i.e.,  $\xi\in
\partial_{cl}f(p),$ which conludes the proof. \hfill
$\diamondsuit$

\begin{remark}\rm If
$(M,g)$ is a Hadamard manifold, then $C_p=\emptyset$ for every $p\in
M$; in this case, the cut-locus subdifferential agrees formally with
the convex subdifferential in the Euclidean setting.
\end{remark}

Let $K\subset M$ be a closed set. Following Ledyaev and Zhu
\cite{LZ-TAMS}, the {\it Fr\'echet-normal cone} and {\it limiting
normal cone} of $K$ at $p\in K$ are the sets $$N_F(p;K)=\partial_F
\delta_K(p)\ \ {\rm and}\ \ N_L(p;K)=\partial_L \delta_K(p),$$
 where $\delta_K$ is the
indicator function of the set $K,$ i.e., $\delta_K(q)=0$ if $q\in K$
and $\delta_K(q)=+\infty$ if $q\notin K$.

The following result - which is one of our key tools to study
Nash-Stampacchia equilibrium points on Riemannian manifolds - it is
know for Hadamard manifolds only, see Li, L\'opez and Mart\'\i
n-M\'arquez \cite{LLMM-1} and it is a simple consequence of the
above theorem.

\begin{corollary}\label{th-normal-cones-tetel} Let $(M,g)$ be a Riemannian
manifold, $K \subset M$ be a closed, geodesic convex set, and $p\in
K$. Then, we have
$$N_F(p;K)=N_L(p;K)=\partial_{cl}\delta_K(p)=\{\xi\in T_p^*M:\langle \xi, \exp_p^{-1}(q)
\rangle_{g}\leq 0\  {for\ all}\  q\in K\}.$$
\end{corollary}

\noindent {\it Proof.} Applying  Theorem
\ref{subdifferentials-theorem-uj} to the indicator function
$f=\delta_K$, we have that
$N_F(p;K)=N_L(p;K)=\partial_{cl}\delta_K(p).$
 It remains to
compute the latter set explicitly. Since $K\subset M\setminus C_p$
(note that the geodesic convexity of $K$ assumes itself that every
two points of $K$ can be joined by a unique geodesic, thus $K\cap
C_p=\emptyset$) and $\delta_K(p)=0$, $\delta_K(q)=+\infty$ for
$q\notin K$, one has that
\begin{eqnarray*}
  \xi\in \partial_{cl}\delta_K(p) &\Leftrightarrow& \delta_K(q)-\delta_K(p)\geq \langle \xi,\exp_p^{-1}(q)\rangle_g\ {\rm for\ all}\ q\in M\setminus C_p \\
   &\Leftrightarrow & 0 \geq \langle \xi,\exp_p^{-1}(q)\rangle_g\ {\rm for\ all}\ q\in
   K,
\end{eqnarray*}\
which ends the proof.  \hfill
$\diamondsuit$\\

Let $U\subset M$ be an open subset of the Riemannian manifold
$(M,g)$. We say that a function $f:U\to \mathbf R$ is {\it locally
Lipschitz at} $p\in U$ if there exist an open neighborhood
$U_p\subset U$ of $p$ and a number $C_p>0$ such that for every
$q_1,q_2\in U_p$, $$|f(q_1)-f(q_2)|\leq C_pd_g(q_1,q_2).$$ The
function $f:U\to \mathbf R$ is {\it locally Lipschitz} on $(U,g)$
 if it is locally Lipschitz at every $p\in U.$

 Fix $p\in U$,  $v\in T_pM$, and let $\tilde U_p\subset U$ be
a totally normal neighborhood of $p.$ If $q\in \tilde U_p$,
following \cite[Section 5]{Azagra-JFA2}, for small values of $|t|$,
we may introduce
$$ 
  \sigma_{q,v}(t)=\exp_q(tw),\ w =d(\exp_q^{-1}\circ
\exp_p)_{\exp_p^{-1}(q)}v.
$$ 
If the function $f:U\to \mathbf R$ is locally Lipschitz on $(U,g)$,
then $$f^0(p;v)=\limsup_{q\to p,\ t\to
0^+}\frac{f(\sigma_{q,v}(t))-f(q)}{t}$$ is called the {\it Clarke
generalized derivative of $f$ at $p\in U$ in direction $v\in T_pM$},
and $$\partial_C f(p)={\rm co}(\partial_L f(p))$$ is the {\it Clarke
subdifferential of $f$ at} $p\in U,$ where 'co' stands for the
convex hull. When $f:U\to \mathbf R$ is a $C^1$ functional at $p\in
U$ then
\begin{equation}\label{minden-differential-hasonlitasa}
    \partial_C f(p)=\partial_L f(p)=\partial_F f(p)=\{df(p)\},
\end{equation}
see \cite[Proposition 4.6]{Azagra-JFA2}. Moreover,  when $(M,g)$ is
the standard Euclidean space, the Clarke subdifferential and the
Clarke generalized gradient agree, see Clarke \cite{Clarke}.

 One can easily prove that the function
$f^0(\cdot;\cdot)$ is upper-semicontinuous on $TU=\cup_{p\in U}T_pM$
and $f^0(p;\cdot)$ is positive homogeneous and subadditive on
$T_pM$, thus convex. In addition, if $U\subset M$ is geodesic convex
and $f:U\to \mathbf R$ is convex, then
\begin{equation}\label{masodik}
f^0(p;v)=\lim_{t\to 0^+}\frac{f(\exp_p(tv))-f(p)}{t},
\end{equation}
 see Claim 5.4 and the first
relation on p. 341 of \cite{Azagra-JFA2}.

\begin{proposition}\label{lipschitz-characterization}{\rm \cite[Corollary
5.3]{LZ-TAMS}} {\it Let $(M,g)$ be a  Riemannian manifold and let
$f:M\to \mathbf R\cup \{+\infty\}$ be a lower semicontinuous
function. Then the following statements are equivalent:
\begin{itemize}
  \item[{\rm (i)}] $f$ is locally Lipschitz at $p\in M;$
  \item[{\rm (ii)}] $\partial_C f$ is bounded in a neighborhood of $p\in
  M;$
  \item[{\rm (iii)}] $\partial_\infty f(p)=\{0\}.$
\end{itemize}}
\end{proposition}

\begin{proposition}\label{prop-konvex-hasno} Let $f,g:M\to \mathbf R\cup\{+\infty\}$ be two
proper,  lower semicontinuous functions. Then, for every $p\in {\rm
dom}(f)\cap {\rm dom}(g)$  with $\partial_{cl} f(p)\neq
\emptyset\neq
\partial_{cl} g(p)$ we have $\partial_{cl} f(p)+\partial_{cl}
g(p)\subset
\partial_{cl}(f+g)(p).$
Moreover, if both functions are convex and $f$ is locally bounded,
the inclusion is  equality.
\end{proposition}

\noindent {\it Proof.} The first part is trivial. For the second
part, $f$ is a locally Lipschitz function (see Azagra, Ferrera, and
L\'opez-Mesas \cite[Proposition 5.2]{Azagra-JFA2}), thus Theorem
\ref{th-normal-cones-tetel} and Propositions \ref{sum-rule} \&
\ref{lipschitz-characterization}  give $
 \partial_{cl}(f+g)(p)  \subset \partial_L(f+g)(p)
   \subset  \partial_L f(p)+\partial_L g(p)
   =\partial_{cl} f(p)+\partial_{cl} g(p).
$ \hfill $\diamondsuit$\\

Let $f:U\to \mathbf R$ be a  locally Lipschitz function and  $p\in
U$. We consider the {\it Clarke 0-subdifferential} of $f$ at $p$ as
\begin{eqnarray*}
\partial_0f(p) &=& \{\xi\in T_p^*M:f^0(p;\exp_p^{-1}(q))\geq \langle \xi,\exp_p^{-1}(q)\rangle_g\ {\rm for\ all}\ q\in U\setminus C_p\} \\
   &=&\{\xi\in T_p^*M:f^0(p;v)\geq \langle \xi,v\rangle_g\ {\rm for\ all}\ v\in
   T_pM\}.
\end{eqnarray*}

\begin{theorem}\label{prop-hasonlitas-subdiff-ujj} Let $(M,g)$ be a Riemannian manifold, $U\subset
M$ be open, $f:U\to \mathbf R$ be a locally Lipschitz function, and
$p\in U$. Then,
 $$\partial_0f(p)=\partial_{cl} (f^0(p;\exp_{p}^{-1}(\cdot)))(p)=\partial_L (f^0(p;\exp_{p}^{-1}(\cdot)))(p)=\partial_C f(p).$$
\end{theorem}
\noindent {\it Proof.} {\bf Step 1.}  $\partial_0f(p)=\partial_{cl}
(f^0(p;\exp_{p}^{-1}(\cdot)))(p).$ It follows  from the definitions.

{\bf Step 2.}   $\partial_0f(p)=\partial_L
(f^0(p;\exp_{p}^{-1}(\cdot)))(p).$\\ The inclusion $"\subset"$
follows from Step 1 and Theorem \ref{subdifferentials-theorem-uj}.
For the converse, we notice that $f^0(p;\exp_{p}^{-1}(\cdot))$ is
locally Lipschitz in a neighborhood of $p$; indeed, $f^0(p;\cdot)$
is convex on $T_pM$ and $\exp_p$ is a local diffeomorphism on a
neighborhood of the origin of $T_pM$.
 Now, let $\xi \in
\partial_L (f^0(p;\exp_{p}^{-1}(\cdot)))(p).$ Then,  $\xi=\lim_k \xi_k$
where $\xi_k\in
\partial_F (f^0(p;\exp_{p}^{-1}(\cdot)))(p_k),\ p_k\to p$.   By Proposition
\ref{subdiff-char}(ii),
for $\psi=\exp_{p}^{-1}:\tilde U_{p}\to T_{p}M$ where $\tilde
U_{p}\subset M$ is a totally normal ball centered at $p$, one has
that {\small
\begin{equation}\label{equation-cut-locus-ujjj} \liminf_{v\to
0}\frac{f^0(p;\exp_{p}^{-1}(p_k)+v)-f^0(p;\exp_{p}^{-1}(p_k))-\langle
\xi_k((d\exp_p)(\exp_p^{-1}(p_k))),v\rangle_g}{\|v\|}\geq
 0.
\end{equation}}
In particular, if $q\in M\setminus C_p$ is fixed arbitrarily and
$v=t\exp_p^{-1}(q)$ for $t>0$ small, the convexity of $f^0(p;\cdot)$
and relation (\ref{equation-cut-locus-ujjj}) yield that
$$f^0(p;\exp_{p}^{-1}(q))\geq \langle
\xi_k((d\exp_p)(\exp_p^{-1}(p_k))),\exp_{p}^{-1}(q)\rangle_g.
$$ Since $\xi=\lim_k \xi_k$, $p_k\to p$ and $d(\exp_p)(0)={\rm id}_{T_pM}$ (see (\ref{exp-identity})), we obtain that $$f^0(p;\exp_{p}^{-1}(q))\geq \langle
\xi,\exp_{p}^{-1}(q)\rangle_g,$$ i.e., $\xi\in \partial_0f(p)$. This
concludes Step 2.

{\bf Step 3.}  $\partial_0f(p)=\partial_C f(p).$\\
First, we prove the inclusion $\partial_0f(p)\subset\partial_C
f(p).$  Here, we follow Borwein and Zhu  \cite[Theorem
5.2.16]{Bor-Zhu}, see also Clarke, Ledyaev, Stern and Wolenski
\cite[Theorem 6.1]{CLSW}. Let $v\in T_pM$ be fixed arbitrarily. The
definition of $f^0(p;v)$ shows that one can choose $t_k\to 0^+$ and
$q_k\to p$ such that
$$f^0(p;v)=\lim_{k\to
\infty}\frac{f(\sigma_{q_k,v}(t_k))-f(q_k)}{t_k}.$$ Fix
$\varepsilon>0$.  For large $k\in \mathbf N$, let $c_k:[0,1]\to M$
be the unique geodesic joining the points $q_k$ and
$\sigma_{q_k,v}(t_k)$, i.e, $c_k(t)=\exp_{q_k}(t
\exp_{q_k}^{-1}(\sigma_{q_k,v}(t_k)))$ and let $U_k=\cup_{t\in
[0,1]} B_g(c_k(t),\varepsilon t_k)$ its $(\varepsilon
t_k)-$neighbor\-hood. Consider also a $C^\infty$ vector field  $V$
on $U_k$ such that $\dot c_k(t)=V(c_k(t))$, $t\in [0,1]$. Now,
applying Proposition \ref{mean-value} with
$r_k=f(c_k(1))-f(c_k(0))-\varepsilon t_k$, one can find
$m_k=m_k(t_k,q_k,v)\in U_k$ and $\xi_k\in \partial_F f(m_k)$ such
that $r_k<\langle \xi_k, V(m_k)\rangle_g$. The latter inequality is
equivalent to
$$\frac{f(\sigma_{q_k,v}(t_k))-f(q_k)}{t_k}<\varepsilon+\left\langle\xi_k,{V(m_k)}/{t_k}\right\rangle_g.$$
Since $f$ is locally Lipschitz, $\partial_F f$ is bounded in a
neighborhood of $p$, see Proposition
\ref{lipschitz-characterization}, thus the sequence $\{\xi_k\}$ is
bounded on $TM$. We can choose a convergent subsequence (still
denoting by $\{\xi_k\}$), and let $\xi_L=\lim_k \xi_k.$ From
construction, $\xi_L\in \partial_L f(p)\subset
\partial_Cf (p).$ Since $m_k\to p$, according to (\ref{exp-identity}), we have that $\lim_{k\to \infty}V(m_k)/t_k=v.$
Thus, letting $k\to \infty$ in the latter inequality, the
arbitrariness of $\varepsilon>0$ yields that
$$f^0(p;v)\leq \langle \xi_L,v\rangle_g.$$
Now, taking into account that $f^0(p;v)=\max
\{\langle\xi,v\rangle_g:\xi\in
\partial_0f(p)\}$, we obtain that
 $$\max
\{\langle\xi,v\rangle_g:\xi\in
\partial_0f(p)\}=f^0(p;v)\leq \langle \xi_L,v\rangle_g\leq \sup\{\langle\xi,v\rangle_g:\xi\in
\partial_Cf(p)\}.$$
H\"ormander's result (see \cite{CLSW}) shows that this inequality in
terms of support functions of convex sets is equivalent to the
inclusion $\partial_0f(p)\subset
\partial_Cf(p)$.

For the converse, it is enough to prove that $\partial_L f(p)\subset
\partial_0f(p)$ since the latter set is convex. Let $\xi\in
\partial_L f(p)$. Then, we have  $\xi=\lim_k \xi_k$ where
$\xi_k\in
\partial_F f(p_k)$ and $ p_k\to p$. A similar argument as in the
proof of Theorem \ref{subdifferentials-theorem-uj} (see relation
(\ref{kell-hatrabbb})) gives that for every $q\in M\setminus C_p$
and $k\in \mathbf N$, we have {\small $$ \liminf_{t\to
0^+}\frac{f(\exp_{p_k}(t \exp_{p_k}^{-1}(q)))-f(p_k)-\langle \xi_k,t
\exp_{p_k}^{-1}(q) \rangle_g}{\|t\exp_{p_k}^{-1}q\|}\geq
 0.
$$}
Since $\|\exp_{p_k}^{-1}q\|=d_g(p_k,q)\geq c_0>0$, by the definition
of the Clarke generalized derivative $f^0$ and the above inequality,
one has that
$$f^0(p_k;\exp_{p_k}^{-1}(q))\geq \langle \xi_k, \exp_{p_k}^{-1}(q)
\rangle_g.$$ The upper semicontinuity of $f^0(\cdot;\cdot)$ and the
fact that $\xi=\lim_k \xi_k$ imply that {\small $$
f^0(p;\exp^{-1}_{p}(q)) \geq  \limsup_k
f^0(p_k;\exp_{p_k}^{-1}(q))\geq \limsup_k\langle \xi_k,
\exp_{p_k}^{-1}(q) \rangle_g= \langle
\xi,\exp_{p}^{-1}(q)\rangle_{g},$$} i.e., $\xi\in
\partial_0f(p)$, which concludes the proof of Step 3.\hfill
$\diamondsuit$\\

\noindent 2.4. {\bf Dynamical systems on manifolds.} In this
subsection we recall the existence of a local solution for a
Cauchy-type problem defined on Riemannian manifolds and its
viability relative to a closed set.


Let $(M,g)$ be a finite-dimensional Riemannian manifold and $G:M\to
{TM}$ be a vector field on $M,$ i.e., $G(p)\in
T_pM$ for every $p\in M$. 
We assume in the sequel that $G:M\to {TM}$ is a $C^{1-0}$ vector
field (i.e., locally Lipschitz); then the dynamical system
$$\leqno{(DS)_G} \ \ \ \ \ \left\{
\begin{array}{lll}
\dot\eta(t)= G(\eta(t)),
\\ \eta(0)={p}_0,
\end{array}
\right. $$ has a unique maximal semiflow $\eta:[0,T)\to M$, see
Chang \cite[p. 15]{Chang}. In particular, $\eta$ is an absolutely
continuous function such that $[0,T)\ni t\mapsto \dot\eta(t)\in
T_{\eta(t)}{M}$ and it verifies $(DS)_G$ for a.e. $t\in [0,T).$


 A set $K\subset M$ is {\it invariant with respect
to the solutions of} $(DS)_G$ if for every initial point $p_0\in K$
the unique maximal semiflow/orbit  $\eta:[0,T)\to M$ of $(DS)_G$
fulfills the property that $\eta(t)\in K$ for every $t\in [0,T).$ We
introduce the  Hamiltonian function as
$$  H_G(p,\xi)=\langle\xi, G(p)\rangle_g, \ (p,\xi)\in
M\times T_p^*M.
$$
Note that $H_G(p,dh(p))<\infty$ for every $p\in M$ and $h\in
C^1(M)$. Therefore, after a suitable adaptation of the results from
Ledyaev and Zhu \cite[Subsection 6.2]{LZ-TAMS} we may state


\begin{proposition}\label{weak-inva-mas} {\it Let $G:M\to
{TM}$ be a $C^{1-0}$ vector field and $K\subset M$ be a non-empty
closed set. The  following statements are equivalent:
\begin{enumerate}
 \item[{\rm (i)}] $K$ is invariant with respect to the solutions of $(DS)_G$;
 \item[{\rm (ii)}] $H_G(p,\xi)\leq 0$ for any $p\in K$ and $\xi\in N_F(p;K)$.
\end{enumerate}}
\end{proposition}

 \section{Comparison of Nash-type equilibria}\label{exist-Nash}

Let $K_1,...,K_n\ (n\geq 2)$ be non-empty sets, corresponding to the
strategies of $n$ players and $f_i:K_1\times...\times K_n
\rightarrow \mathbf R$ $(i\in \{1,... ,n\})$ be the payoff
functions, respectively. Throughout the paper, the following
notations/conventions are used:
\begin{itemize}
  \item[$\bullet$] ${\bf K}=K_1\times ... \times K_n;$ ${ \bf f}=(f_1,...,f_n);$ ${ \bf (f,K)}=(f_1,...,f_n;K_1,...,K_n);$
  \item[$\bullet$] ${\bf p}=(p_1,...,p_n)$;
  \item[$\bullet$] ${\bf p}_{-i}$ is a strategy profile of all players except for player $i;$\\ $(q_i,{\bf  p}_{-i})=(p_1,...,p_{i-1},q_i,p_{i+1},...,p_n);$ in
  particular, $(p_i,{\bf  p}_{-i})={\bf  p};$
  \item[$\bullet$] ${\bf K}_{-i}$ is the strategy set profile of all players except for player $i;$\\
  $(U_i,{\bf K}_{-i})=K_1\times...\times K_{i-1} \times U_i\times K_{i+1}\times...\times K_n$
for some $U_i\supset K_i.$
\end{itemize}

\noindent
\begin{definition}\label{Nash-equil}
The set of {\it Nash equilibrium points for ${ \bf (f,K)}$} is
 $${\mathcal S}_{NE}{ \bf (f,K)}=\big\{{\bf p}\in {\bf K}:f_i(q_i,{\bf
p}_{-i})\geq f_i({\bf p})\ \ {\rm for\ all}\ q_i\in K_i,\ i\in
\{1,... ,n\}\big\}.$$
\end{definition}
 \noindent The main result of the paper
\cite{Kristaly-PAMS-2009} states that in a quite general framework
the set of Nash equilibrium points for ${ \bf (f,K)}$ is not empty.
More precisely, we have

\begin{proposition}\label{main1-uj}{\rm \cite{Kristaly-PAMS-2009}}
{\it Let $(M_i,g_i)$ be finite-dimensional Riemannian manifolds;
$K_i\subset M_i$ be non-empty, compact, geodesic convex sets; and
$f_i:{\bf K}\to \mathbf R$ be continuous functions such that $K_i\ni
q_i\mapsto f_i(q_i,{\bf p}_{-i})$ is convex on $K_i$ for every ${\bf
p}_{-i}\in {\bf K}_{-i},$ $i\in \{1,... ,n\}$. Then there exists at
least one Nash equilibrium point for ${\bf (f,K)}$, i.e., ${\mathcal
S}_{NE}{ \bf (f,K)}\neq \emptyset$.}
\end{proposition}

Similarly to Proposition \ref{main1-uj}, let us assume that for
every $i\in \{1,... ,n\}$, one can find a finite-dimensional
Riemannian manifold $(M_i,g_i)$ such that the strategy set $K_i$ is
closed and geodesic convex in $(M_i,g_i)$.  Let ${\bf M}=
M_1\times...\times M_n$ be the product manifold with its standard
 Riemannian product metric
\begin{equation}\label{prod-metric}
 {\bf g(V,W)}=\sum_{i=1}^n g_i(V_i,W_i)
\end{equation}
for every ${\bf V}=(V_1,...,V_n), {\bf W}=(W_1,...,W_n)\in
T_{p_1}M_1\times ...\times T_{p_n}M_n=T_{\bf p}{\bf M}$. Let ${\bf
U}=U_1\times ...\times U_n\subset {\bf M}$ be an open set such that
${\bf K}\subset {\bf U}$; we always mean that $U_i$ inherits the
Riemannian structure of $(M_i,g_i)$. Let
\begin{eqnarray*}
 \mathcal L_{({\bf K, U, M})}=\big\{{\bf f}=(f_1,...,f_n)\in C^0({\bf K}, \mathbf
 R^n)
& : & f_i:(U_i,{\bf K}_{-i} )\to \mathbf R\ \mbox{is continuous
and}\\&& f_i(\cdot,{\bf p}_{-i})\ \mbox{is locally Lipschitz}\
\mbox{on} \ (U_i,g_i)\\ && \mbox{for all}\ {\bf p}_{-i}\in {\bf
K}_{-i}, \ i\in \{1,... ,n\}\big\}.
\end{eqnarray*}
The next notion has been introduced in \cite{Kristaly-PAMS-2009}.

\begin{definition} Let ${\bf f}\in \mathcal L_{({\bf K, U,M})}.$ The set of {\it Nash-Clarke equilibrium points for} ${ \bf (f,K)}$
is $${\mathcal S}_{NC}({\bf f,K})=\big\{{\bf p}\in {\bf K}:
f_i^0({\bf p};\exp^{-1}_{p_i}(q_i)) \geq 0\ \ {\rm for\ all}\ q_i\in
K_i,\ i\in \{1,... ,n\}\big\}.$$
\end{definition}
\noindent Here, $f_i^0({\bf p};\exp^{-1}_{p_i}(q_i))$ denotes the
Clarke generalized derivative of $f_i(\cdot,{\bf p}_{-i})$ at point
$p_i\in K_i$ in direction $\exp^{-1}_{p_i}(q_i)\in T_{p_i}M_i.$ More
precisely,
\begin{equation}\label{iranymenti-deri-clarke}
  f_i^0({\bf p};\exp^{-1}_{p_i}(q_i))=\limsup_{q\to p_i,q\in U_i,\
t\to 0^+}\frac{f_i(\sigma_{q,\exp^{-1}_{p_i}(q_i)}(t),{\bf p
}_{-i})-f_i(q,{\bf p }_{-i})}{t},
\end{equation}
where
 $ \sigma_{q,v}(t)=\exp_q(tw),$ and $ w =d(\exp_q^{-1}\circ
\exp_{p_i})_{\exp_{p_i}^{-1}(q)}v$ for $v\in T_{p_i}M_i$, and $t>0$
is small enough.
By exploiting a minimax result of McClendon \cite{McClendon}, the
following existence result is available concerning the Nash-Clarke
points for ${\bf (f,K)}$.

\begin{proposition}\label{main2}{\rm \cite{{Kristaly-PAMS-2009}}}
{\it Let $(M_i,g_i)$ be complete finite-dimensional Riemannian
manifolds; $K_i\subset M_i$ be non-empty, compact, geodesic convex
sets; and ${\bf f}\in \mathcal L_{({\bf K, U,M})}$  such that for
every ${\bf p}\in {\bf K}$, $i\in \{1,... ,n\}$, $K_i\ni q_i\mapsto
f_i^0({\bf p};\exp_{p_i}^{-1}(q_i))$ is convex and $ f_i^0$ is upper
semicontinuous on its domain of definition. Then ${\mathcal
S}_{NC}({\bf f,K})\neq \emptyset$.}
\end{proposition}

\begin{remark}\rm Although Proposition \ref{main2} gives a possible
approach to locate Nash equilibria on Riemannian
manifolds, 
its applicability is quite reduced. Indeed, ${f_i}^0({\bf
p};\exp_{p_i}^{-1}(\cdot))$ has no convexity property in general,
unless we are in the Euclidean setting or the set $K_i$ is a
geodesic segment, see \cite{Kristaly-PAMS-2009}. For instance, if
$\mathbf H^2$ is the standard Poincar\'e upper-plane with the metric
$g_{\mathbf H}=(\frac{\delta_{ij}}{y^2})$ and we consider the
function $f:\mathbf H^2\times \mathbf R\to \mathbf R$,
$f((x,y),r)=rx$ and the geodesic segment $\gamma(t)=(1,e^t)$ in
$\mathbf H^2$, $t\in [0,1]$, the function {\small $$t\mapsto
f_1^0(((2,1),r);\exp_{(2,1)}^{-1}(\gamma(t)))=r\left(e^{2t}\frac{\sinh
2}{2}+e^t\cosh 1\sqrt{e^{2t}(\cosh 1)^2-1} \right)^{-1}$$} is  not
convex.
\end{remark}

The limited applicability of Proposition \ref{main2} comes from the
involved form of the set ${\mathcal S}_{NC}{ \bf (f,K)}$ which
motivates the introduction and study of the following concept which
plays the central role in the present paper.

\begin{definition}\label{Nash-Stampacchia}
Let ${\bf f}\in \mathcal L_{({\bf K, U, M})}.$ The set of {\it
Nash-Stampacchia equilibrium points for ${ \bf (f,K)}$} is
\begin{eqnarray*}
{\mathcal S}_{NS}{ \bf (f,K)}=\big\{{\bf p}\in {\bf K}&:&\exists
\xi_C^i\in
\partial_C^i f_i({\bf p}) \ \mbox{such that}\ \langle \xi_C^i,
\exp_{p_i}^{-1}(q_i)\rangle_{g_i} \geq 0,\\ &&  {\rm for\ all}\
q_i\in K_i,\ i\in \{1,... ,n\}\big\}.
\end{eqnarray*}
\end{definition}
\noindent Here, $\partial_C^i f_i({\bf p})$ denotes the Clarke
subdifferential of the function $f_i(\cdot,{\bf p}_{-i})$ at point
$p_i\in K_i,$ i.e., $\partial_C f_i(\cdot,{\bf p}_{-i})(p_i)={\rm
co}(\partial_L f_i(\cdot,{\bf p}_{-i})(p_i))$.\\

Our first aim is to compare the three Nash-type points introduced in
Definitions \ref{Nash-equil}-\ref{Nash-Stampacchia}. Before to do
that, we introduce another two classes of functions. If
${U_i}\subset M_i$ is geodesic convex for every $i\in \{1,... ,n\}$,
we may define
\begin{eqnarray*}
 \mathcal K_{({\bf K, U,M})}=\big\{{\bf f}\in C^0({\bf K},\mathbf R^n) & : &
f_i:(U_i,{\bf K}_{-i})\to \mathbf R\ \mbox{is continuous and}\
f_i(\cdot,{\bf p}_{-i})\ \mbox{is} \\ && \mbox{convex}\ \mbox{on} \
(U_i,g_i)\ \mbox{for all}\ {\bf p}_{-i}\in {\bf K}_{-i}, \ i\in
\{1,... ,n\}\big\},
\end{eqnarray*}
and
\begin{eqnarray*}
 \mathcal C_{({\bf K, U,M})}=\big\{{\bf f}\in C^0({\bf K},\mathbf R^n) & : &
f_i:(U_i,{\bf K}_{-i})\to \mathbf R\ \mbox{is continuous and}\
f_i(\cdot,{\bf p}_{-i})\ \mbox{is of} \\ && \mbox{class}\ C^1\
\mbox{on} \ (U_i,g_i)\ \mbox{for all}\ {\bf p}_{-i}\in {\bf K}_{-i},
\ i\in \{1,... ,n\}\big\}.
\end{eqnarray*}

\begin{remark}\label{remark-azagra-lok-lip}\rm  Due to Azagra, Ferrera and L\'opez-Mesas \cite[Proposition
5.2]{Azagra-JFA2}, one has  that $\mathcal K_{({\bf K, U,M})}\subset
\mathcal L_{({\bf K, U,M})}$. Moreover, it is clear that $\mathcal
C_{({\bf K, U,M})}\subset \mathcal L_{({\bf K, U,M})}$.

\end{remark}
The main result of this section reads as follows.

\begin{theorem}\it \label{equi-kritikus} Let $(M_i,g_i)$ be finite-dimensional
Riemannian manifolds; $K_i\subset M_i$ be non-empty, closed,
geodesic convex sets;  $U_i\subset M_i$ be open sets containing
$K_i$; and $f_i:{\bf K}\to \mathbf R$ be some functions, $i\in
\{1,... ,n\}$. Then, we have
 \begin{enumerate}
 \item[{\rm (i)}] ${\mathcal S}_{NE}{ \bf (f,K)}\subset {\mathcal S}_{NS}{ \bf (f,K)}={\mathcal S}_{NC}{ \bf (f,K)}$
 whenever ${\bf f}\in \mathcal L_{({\bf K, U, M})};$
  \item[{\rm (ii)}]  ${\mathcal S}_{NE}{ \bf (f,K)}= {\mathcal S}_{NS}{ \bf (f,K)}= {\mathcal S}_{NC}{ \bf (f,K)}$
 whenever ${\bf f}\in \mathcal K_{({\bf K, U, M})};$
 \end{enumerate}
\end{theorem}

\noindent {\it Proof.} (i) First, we prove that ${\mathcal S}_{NE}{
\bf (f,K)}\subset {\mathcal S}_{NS}{ \bf (f,K)}$. Indeed, we have\\
${\bf p}\in {\mathcal S}_{NE}{ \bf (f,K)}\Leftrightarrow$
\begin{eqnarray*}
  &\Leftrightarrow& f_i(q_i,{\bf
p}_{-i})\geq f_i({\bf p})\ \ {\rm for\ all}\ q_i\in K_i,\ i\in
\{1,... ,n\}\\
  &\Leftrightarrow & 0\in \partial_{cl}(f_i(\cdot,{\bf p}_{-i})+\delta_{K_i})(p_i),\ i\in\{1,...,n\} \\
  &\Rightarrow & 0\in \partial_L(f_i(\cdot,{\bf p}_{-i})+\delta_{K_i})(p_i),\ i\in\{1,...,n\} \ \ \ \ ({\rm cf.\ Theorem\ \ref{subdifferentials-theorem-uj}})\\
  &\Rightarrow & 0\in \partial_Lf_i(\cdot,{\bf p}_{-i})(p_i)+\partial_L \delta_{K_i}(p_i),\ i\in\{1,...,n\} \ \ \ \ ({\rm cf.\ Propositions\ \ref{sum-rule}\ \&\ \ref{lipschitz-characterization}})\\
  &\Rightarrow & 0\in \partial_Cf_i(\cdot,{\bf p}_{-i})(p_i)+\partial_L \delta_{K_i}(p_i),\ i\in\{1,...,n\} \\
 &\Leftrightarrow & 0\in \partial_C^i {f_i}({\bf p})+N_L(p_i;K_i),\ i\in\{1,...,n\}\\
 &\Leftrightarrow& \exists
\xi_C^i\in \partial_C^i {f_i}({\bf p})\ \mbox{such that}\
\nonumber\langle\xi_C^i,\exp_{p_i}^{-1}({q_i})\rangle_{{g_i}}\geq 0\ \ {\rm for\ all}\ {q_i}\in {K_i},i\in\{1,...,n\}\\
&&\ \ \ \
 \ \ \ \ \ \ \ \ \ \ \ \ \ \ \ \ \ \ \ \ \
 \ \ \ \ \ \ \ \ \ \ \ \ \ \ \ \ \ \ \ \ \
\ \ \ \
 \ \ \ \ \ \ \ \ \ \ \ \ \ \ \ \ \  \ \ \ \ \ \ \ \ ({\rm cf.\ Corollary\ \ref{th-normal-cones-tetel}})\\
  &\Leftrightarrow & {\bf p}\in {\mathcal S}_{NS}{ \bf (f,K)}.
\end{eqnarray*}
Now, we prove ${\mathcal S}_{NS}{ \bf (f,K)}\subset {\mathcal
S}_{NC}{ \bf (f,K)}$; more precisely, we have\\ $ {\bf p}\in
{\mathcal S}_{NS}{ \bf (f,K)}\Leftrightarrow$
\begin{eqnarray*}
 &\Leftrightarrow & 0\in \partial_C^i {f_i}({\bf p})+N_L(p_i;K_i),\ i\in\{1,...,n\}\\
  &\Leftrightarrow & 0\in \partial_C^i {f_i}({\bf p})+\partial_{cl}\delta_{K_i}(p_i),\ i\in\{1,...,n\}\ \ \ \ ({\rm cf.\ Corollary\ \ref{th-normal-cones-tetel}}) \\
  &\Leftrightarrow & 0\in \partial_{cl} ({f_i}^0({\bf p};\exp_{p_i}^{-1}(\cdot)))(p_i)+\partial_{cl}\delta_{K_i}(p_i),\ i\in\{1,...,n\} \ \ \ \ ({\rm cf.\ Theorem\ \ref{prop-hasonlitas-subdiff-ujj}}) \\
  &\Rightarrow & 0\in \partial_{cl} ({f_i}^0({\bf p};\exp_{p_i}^{-1}(\cdot))+\delta_{K_i})(p_i),\ i\in\{1,...,n\} \ \ \ \ ({\rm cf.\ Proposition\ \ref{prop-konvex-hasno} 
  })\\
  &\Leftrightarrow & f_i^0({\bf p};\exp_{p_i}^{-1}(q_i))\geq 0\ \ {\rm for\ all}\  q_i\in K_i,\ i\in\{1,...,n\} \\
  &\Leftrightarrow & {\bf p}\in {\mathcal S}_{NC}{ \bf (f,K)}.
\end{eqnarray*}
In order to prove  ${\mathcal S}_{NC}{ \bf (f,K)}\subset {\mathcal
S}_{NS}{ \bf (f,K)}$, we recall that  ${f_i}^0({\bf
p};\exp_{p_i}^{-1}(\cdot))$ is locally Lipschitz in a neighborhood
of $p_i$. Thus, we have
 \\
$ {\bf p}\in {\mathcal S}_{NC}{ \bf
(f,K)}\Leftrightarrow$\begin{eqnarray*}
 &\Leftrightarrow & 0\in \partial_{cl} ({f_i}^0({\bf p};\exp_{p_i}^{-1}(\cdot))+\delta_{K_i})(p_i),\
 i\in\{1,...,n\}\\
  &\Rightarrow & 0\in \partial_L ({f_i}^0({\bf p};\exp_{p_i}^{-1}(\cdot))+\delta_{K_i})(p_i),\ i\in\{1,...,n\} \ \ \ ({\rm cf.\ Theorem\ 2.1}) \\
  &\Rightarrow & 0\in \partial_L ({f_i}^0({\bf p};\exp_{p_i}^{-1}(\cdot)))(p_i)+\partial_L\delta_{K_i}(p_i),\ i\in\{1,...,n\} \ \
  ({\rm cf.\ Propositions\ \ref{sum-rule}\ \&\ \ref{lipschitz-characterization}})\\
  &\Leftrightarrow & 0\in \partial_C ({f_i}(\cdot, {\bf p_{-i}}))(p_i)+\partial_L\delta_{K_i}(p_i),\ i\in\{1,...,n\} \ \ \ \ ({\rm cf.\
  Theorem\ \ref{prop-hasonlitas-subdiff-ujj}
  })\\
  &\Leftrightarrow & 0\in \partial_C^i ({f_i}( {\bf p}))+N_L(p_i;K_i),\
  i\in\{1,...,n\}\\
  &\Leftrightarrow & {\bf p}\in {\mathcal S}_{NS}{ \bf (f,K)}.
\end{eqnarray*}

(ii) 
Due to (i) and Remark \ref{remark-azagra-lok-lip}, it is enough to
prove that ${\mathcal S}_{NC}{ \bf (f,K)}\subset {\mathcal S}_{NE}{
\bf (f,K)}$. Let ${\bf p}\in {\mathcal S}_{NC}{ \bf (f,K)}$, i.e.,
for every $i\in \{1,... ,n\}$ and $q_i\in K_i$,
\begin{equation}\label{sssms}
  f_i^0({\bf p};\exp^{-1}_{p_i}(q_i)) \geq 0.
\end{equation}
Fix  $i\in \{1,... ,n\}$ and $q_i\in K_i$ arbitrary. Since $
f_i(\cdot,{\bf p}_{-i})$ is convex on $(U_i,g_i)$, on account of
(\ref{masodik}), we have
\begin{equation}\label{g(t)-Nash}
   f_i^0({\bf p};\exp_{p_i}^{-1}(q_i))=\lim_{t\to
0^+}\frac{f_i(\exp_{p_i}(t\exp_{p_i}^{-1}(q_i)),{\bf
p}_{-i})-f_i({\bf p})}{t}.
\end{equation}
Note that the function
$$R(t)=\frac{f_i(\exp_{p_i}(t\exp_{p_i}^{-1}(q_i)),{\bf
p}_{-i})-f_i({\bf p})}{t}$$ is well-defined on the whole interval
$(0,1]$; indeed, $t\mapsto \exp_{p_i}(t\exp_{p_i}^{-1}(q_i))$ is the
minimal geodesic joining the points $p_i\in K_i$ and $q_i\in K_i$
which belongs to $K_i\subset U_i.$ Moreover, it is well-known  that
$t\mapsto R(t)$ is non-decreasing on $(0,1]$. Consequently,
$$f_i(q_i,{\bf p}_{-i})-f_i({\bf p})=f_i(\exp_{p_i}(\exp_{p_i}^{-1}(q_i)),{\bf
p}_{-i})-f_i({\bf p})=R(1)\geq\lim_{t\to 0^+}R(t).$$ Now,
(\ref{sssms}) and (\ref{g(t)-Nash}) give that $\lim_{t\to
0^+}R(t)\geq 0,$ which concludes the proof.
 \hfill $\diamondsuit$

\begin{remark}\rm (a) As we can see, the key tool in the proof of ${\mathcal S}_{NS}{ \bf (f,K)}={\mathcal S}_{NC}{ \bf (f,K)}$ is
the locally Lipschitz property of the function ${f_i}^0({\bf
p};\exp_{p_i}^{-1}(\cdot))$  near $p_i$.

(b) In \cite{Kristaly-PAMS-2009} we considered the sets ${\mathcal
S}_{NE}({\bf f,K})$ and ${\mathcal S}_{NC}({\bf f,K})$. Note however
that the Nash-Stampacchia concept is more appropriate to find Nash
equilibrium points in  general contexts, see also the applications
in \S \ref{utol-sect} for both compact and non-compact cases.
Moreover, via ${\mathcal S}_{NS}({\bf f,K})$ we realize that the
optimal geometrical framework to develop this study is the class of
Hadamard manifolds. In the next  sections we develop this approach.
\end{remark}

\section{Nash-Stampacchia equilibria on Hadamard manifolds: characterization, existence and stability}\label{section-main}

Let $(M_i,g_i)$ be finite-dimensional Hadamard manifolds, $i\in
\{1,...,n\}$. Standard arguments show that $({\bf M,g})$ is also a
Hadamard manifold, see Ballmann \cite[Example 4, p.147]{Ballmann}
and O'Neill \cite[Lemma 40, p. 209]{ONeill}. Moreover, on account of
the characterization of (warped) product geodesics, see O'Neill
\cite[Proposition 38, p. 208]{ONeill}, if $\exp_{\bf p}$ denotes the
usual exponential map on $({\bf M,g})$ at ${\bf p}\in {\bf M},$ then
for every $ {\bf V}=(V_1,...,V_n)\in T_{\bf p}{\bf M}$, we have
$$\exp_{\bf p}({\bf V})=(\exp_{p_1}(V_1),...,\exp_{p_n}(V_n)).$$ We
consider that $K_i\subset M_i$ are non-empty, closed, geodesic
convex sets and $U_i\subset M_i$ are open sets containing $K_i$,
$i\in \{1,...,n\}$.

Let ${\bf f}\in \mathcal L_{({\bf K, U, M})}.$ The {\it diagonal
Clarke subdifferential of} ${\bf f}=(f_1,...,f_n)$ at ${\bf p}\in
{\bf K}$ is
 $$\partial_C^\Delta {\bf f}({\bf p})=(\partial_C^1
f_1({\bf p}),...,\partial_C^n f_n({\bf p})).$$
From the definition of the metric ${\bf g}$, for every ${\bf p}\in
{\bf K}$ and ${\bf q}\in {\bf M}$ it turns out that
\begin{equation}\label{metrika-osszef}
  \langle\xi_C^\Delta,\exp_{\bf p}^{-1}({\bf q})\rangle_{{\bf g}}=\sum_{i=1}^n
\langle\xi_C^i,\exp_{p_i}^{-1}(q_i)\rangle_{g_i},\ \
\xi_C^\Delta=(\xi_C^1,...,\xi_C^n)\in \partial_C^\Delta {\bf f}({\bf
p}).
\end{equation}

\vspace{0.3cm} \noindent 4.1. {\bf Nash-Stampacchia equilibrium
points versus fixed points of $A_\alpha^{\bf f}$.} For each
$\alpha>0$ and ${\bf f}\in \mathcal L_{({\bf K, U, M})}$, we define
the set-valued map $A_\alpha^{\bf f}:{\bf K}\to 2^{\bf K}$ by
$$A_\alpha^{\bf f}({\bf p})=P_{\bf K}(\exp_{\bf p}(-\alpha\partial_C^\Delta
{\bf f}({\bf p}))), \ \ {\bf p\in K}.$$
 Note that for each ${\bf p\in K}$, the set $A_\alpha^{\bf f}({\bf p})$ is non-empty and compact.   The following result plays a crucial role in our further
 investigations.
\begin{theorem}\label{theorem-equivalence} Let $(M_i,g_i)$ be finite-dimensional
Hadamard manifolds; $K_i\subset M_i$ be non-empty, closed, geodesic
convex sets;  $U_i\subset M_i$ be open sets containing $K_i$, $i\in
\{1,... ,n\}$; and ${\bf f}\in \mathcal L_{({\bf K, U, M})}$. Then
the following statements are equivalent:
\begin{itemize}
  \item[{\rm (i)}] ${\bf p}\in {\mathcal S}_{NS}({\bf f,K});$
  \item[{\rm (ii)}] ${\bf p}\in A_\alpha^{\bf f}({\bf p})$ for all $\alpha>0;$
  \item[{\rm (iii)}] ${\bf p}\in A_\alpha^{\bf f}({\bf p})$ for some $\alpha>0$.
\end{itemize}
\end{theorem}
{\it Proof.} In view of relation (\ref{metrika-osszef}) and the
identification between $T_{\bf p}{\bf M}$ and $T_{\bf p}^*{\bf M}$,
see (\ref{dual-nem-dual}), we have that
\begin{eqnarray}\label{ekviv-1}
{\bf p}\in {\mathcal S}_{NS}({\bf f,K}) &\Leftrightarrow &\exists
\xi_C^\Delta=(\xi_C^1,...,\xi_C^n)\in \partial_C^\Delta {\bf f}({\bf
p})\ \mbox{such that}\\ && \nonumber\langle\xi_C^\Delta,\exp_{\bf
p}^{-1}({\bf
q})\rangle_{{\bf g}}\geq 0\ \ {\rm for\ all}\ {\bf q}\in {\bf K}\\
\nonumber  &\Leftrightarrow &\exists
\xi_C^\Delta=(\xi_C^1,...,\xi_C^n)\in
\partial_C^\Delta {\bf f}({\bf p})\ \mbox{such that}\\ &&
\nonumber {\bf g}(-\alpha\xi_C^\Delta,\exp_{\bf p}^{-1}({\bf
q}))\leq 0\ \ {\rm for\ all}\ {\bf q}\in {\bf K}\ {\rm and}\\&& \
{\rm for\ all/some}\ \alpha>0.\nonumber 
\end{eqnarray}
On the other hand, let $\gamma,\sigma:[0,1]\to {\bf M}$
 be the unique minimal geodesics defined
by $\gamma(t)=\exp_{\bf p}(-t\alpha\xi_C^\Delta)$ and
$\sigma(t)=\exp_{\bf p}(t\exp_{\bf p}^{-1}({\bf q}))$ for any fixed
$\alpha>0$ and ${\bf q}\in {\bf K}$. Since ${\bf K}$ is geodesic
convex in $({\bf M,g})$, then Im$\sigma\subset {\bf K}$ and
\begin{equation}\label{masodik-2}
  {\bf g}(\dot\gamma(0),\dot\sigma(0))={\bf g}(-\alpha\xi_C^\Delta,\exp_{\bf p}^{-1}({\bf q})).
\end{equation}
Taking into account relation (\ref{masodik-2}) and Proposition
\ref{egyik-irany} (i), i.e., the validity of the obtuse-angle
property on the Hadamard manifold $({\bf M,g})$, (\ref{ekviv-1}) is
equivalent to
$$ {\bf p}=\gamma(0)=  P_{\bf K}(\gamma(1))=P_{\bf K}(\exp_{\bf
 p}(-\alpha\xi_C^\Delta)),$$
which is nothing but ${\bf p} \in A_\alpha^{\bf f}({\bf p})$.
 \hfill $\diamondsuit$

\begin{remark}\rm
Note that the implications {\rm (ii)}$\Rightarrow${\rm (i)} and {\rm
(iii)}$\Rightarrow${\rm (i)} hold for arbitrarily Riemannian
manifolds, 
see Remark \ref{remark-obtuse-kicsi} (b). These implications are
enough to find Nash-Stampacchia equilibrium points for $({\bf f,K})$
via fixed points of the map $A_\alpha^{\bf f}$. However, in the
sequel we exploit further aspects of the Hadamard manifolds as
non-expansiveness of the projection operator of geodesic convex sets
and a Rauch-type comparison property. Moreover, in the spirit of
Nash's original idea that Nash equilibria appear exactly as fixed
points of a specific map, Theorem \ref{theorem-equivalence} provides
a full characterization of Nash-Stampacchia equilibrium points for
$({\bf f,K})$ via the fixed points of the set-valued map
$A_\alpha^{\bf f}$ when $(M_i,g_i)$ are Hadamard manifolds.
\end{remark}

In the sequel, two cases will be considered  to guarantee
Nash-Stampacchia equilibrium points for $({\bf f,K})$, depending on
the compactness of the strategy sets $K_i$.

\vspace{0.3cm}

\noindent 4.2. {\bf Nash-Stampacchia equilibrium points; compact
case.} Our first result guarantees the existence of a
Nash-Stampacchia equilibrium point for $({\bf f,K})$ whenever the
sets $K_i$ are compact; 
the proof is based on Begle's fixed point theorem for set-valued
maps. More precisely, we have

\begin{theorem}\label{begle}
Let $(M_i,g_i)$ be finite-dimensional Hadamard manifolds;
$K_i\subset M_i$ be non-empty, compact, geodesic convex sets; and
$U_i\subset M_i$ be open sets containing $K_i$, $i\in \{1,... ,n\}$.
Assume that ${\bf f}\in \mathcal L_{({\bf K,U,M})}$ and ${\bf K}\ni
{\bf p}\mapsto\partial_C^\Delta {\bf f}({\bf p})$ is upper
semicontinuous.
 Then there exists at
least one Nash-Stampacchia equilibrium point for ${\bf (f,K)}$,
i.e., ${\mathcal S}_{NS}{ \bf (f,K)}\neq \emptyset$.
\end{theorem}

\noindent {\it Proof.} Fix $\alpha>0$ arbitrary. We prove that the
set-valued map $A_\alpha^{\bf f}$ has closed graph. Let $({\bf p
},{\bf q})\in {\bf K}\times {\bf K}$ and the sequences $\{{\bf p
}_k\}, \{{\bf q}_k\}\subset {\bf K}$ such that ${\bf q}_k\in
A_\alpha^{\bf f}({\bf p }_k)$ and $({\bf p }_k, {\bf q}_k)\to ({\bf
p },{\bf q})$ as $k\to \infty.$ Then, for every $k\in \mathbf N$,
there exists $\xi_{C,k}^\Delta\in
\partial_C^\Delta {\bf f}({\bf p}_k)$ such that ${\bf q}_k=P_{\bf K}(\exp_{{\bf
p}_k}(-\alpha\xi_{C,k}^\Delta))$. On account of Proposition
\ref{lipschitz-characterization} (i)$\Leftrightarrow$(ii), the
sequence $\{\xi_{C,k}^\Delta\}$ is bounded on the cotangent bundle
$T^*{\bf M}$. Using the identification between elements of the
tangent and cotangent fibers, up to a subsequence, we may assume
that $\{\xi_{C,k}^\Delta\}$ converges to an element
$\xi_{C}^\Delta\in T_{\bf p}^*{\bf M}$. Since the set-valued map
$\partial_C^\Delta {\bf f}$ is upper semicontinuous on ${\bf K}$ and
${\bf p}_k\to {\bf p}$ as $k\to \infty$, we  have that
$\xi_{C}^\Delta\in \partial_C^\Delta {\bf f}({\bf p})$.  The
non-expansiveness of $P_{\bf K}$ (see Proposition \ref{egyik-irany}
(ii)) gives that
\begin{eqnarray*}
{\bf d_g}({\bf q},P_{\bf K}(\exp_{\bf
p}(-\alpha\xi_{C}^\Delta)))&\leq& {\bf d_g}({\bf q},{\bf q}_k)+{\bf
d_g}({\bf q}_k,P_{\bf K}(\exp_{\bf
p}(-\alpha\xi_{C}^\Delta)))\\&=&{\bf d_g}({\bf q},{\bf q}_k)+{\bf
d_g}(P_{\bf K}(\exp_{{\bf p}_k}(-\alpha\xi_{C,k}^\Delta)),P_{\bf
K}(\exp_{\bf p}(-\alpha\xi_{C}^\Delta)))\\&\leq&{\bf d_g}({\bf
q},{\bf q}_k)+{\bf d_g}(\exp_{{\bf
p}_k}(-\alpha\xi_{C,k}^\Delta),\exp_{\bf p}(-\alpha\xi_{C}^\Delta))
\end{eqnarray*}
Letting $k\to \infty$, both terms in the last expression tend to
zero. Indeed, the former follows from the fact that ${\bf q}_k\to
{\bf q}$ as $k\to \infty$, while the latter is a simple consequence
of the local behaviour of the exponential map. Thus, $${\bf
q}=P_{\bf K}(\exp_{\bf p}(-\alpha\xi_{C}^\Delta))\in P_{\bf
K}(\exp_{\bf p}(-\alpha\partial_C^\Delta {\bf f}({\bf
p})))=A_\alpha^{\bf f}({\bf p}),$$ i.e., the graph of $A_\alpha^{\bf
f}$ is closed.

By definition, for each ${\bf p}\in {\bf K}$ the set
$\partial_C^\Delta {\bf f}({\bf p})$ is convex, so contractible.
Since both $P_{\bf K}$ and the exponential map are continuous,
$A_\alpha^{\bf f}({\bf p})$ is contractible as well for each ${\bf
p}\in {\bf K}$, so acyclic (see \cite{McClendon}).

Now, we are in position to apply Begle's fixed point theorem, see
for instance McClendon \cite[Proposition 1.1]{McClendon}.
Consequently, there exists ${\bf p}\in {\bf K }$ such that ${\bf
p}\in A_\alpha^{\bf f}({\bf p})$. On account of Theorem
\ref{theorem-equivalence}, ${\bf p}\in {\mathcal S}_{NS}({\bf
f,K}).$ \hfill $\diamondsuit$



\begin{remark}\rm
(a) Since ${\bf f}\in \mathcal L_{({\bf K,U,M})}$ in Theorem
\ref{begle}, the partial Clarke gradients $q\mapsto\partial_C
f_i(\cdot,{\bf p}_{-i})(q)$ are upper semicontinuous, $i\in
\{1,...,n\}$. However, in general, the diagonal Clarke
subdifferential $\partial_C^\Delta {\bf f}(\cdot)$ does not inherit
this regularity property.

(b) Two unusual applications to Theorem \ref{begle} will be given in
Examples \ref{example-1} and \ref{example-2}; the first on the
Poincar\'e disc, the second on the manifold of positive definite,
symmetric matrices.
\end{remark}
\vspace{0.3cm}

\noindent 4.3. {\bf Nash-Stampacchia equilibrium points; non-compact
case.} In the sequel, we are focusing to the location of
Nash-Stampacchia equilibrium points for $(\bf f,K)$ in the case when
$K_i$ are {\it not} necessarily compact on the Hadamard manifolds
$(M_i,g_i)$. Simple examples show that even the
$C^\infty-$smoothness of the payoff functions are not enough to
guarantee the existence of Nash(-Stampacchia) equilibria. Indeed, if
$f_1,f_2:\mathbf R^2\to \mathbf R$ are defined as
$f_1(x,y)=f_2(x,y)=e^{-x-y}$, and $K_1=K_2=[0,\infty)$, then
${\mathcal S}_{NS}({\bf f,K})={\mathcal S}_{NE}({\bf
f,K})=\emptyset.$ Therefore, in order to prove existence/location of
Nash(-Stampacchia) equilibria on not necessarily compact strategy
sets, one needs to require more specific assumptions on ${\bf
f}=(f_1,...,f_n)$. Two such possible ways are described in the
sequel.

The first existence result is based on a suitable coercivity
assumption and Theorem \ref{begle}. For
a fixed ${\bf p}_0\in {\bf K}$, we introduce the hypothesis:\\
\\
$(H_{{\bf p}_0})\ \ {\rm There\ exists}\ \xi_C^0\in
\partial_C^\Delta {\bf f}({\bf p}_0)\ {\rm such\ that}$
{\small $$L_{{\bf p}_0}=\limsup_{{\bf d_g}({\bf p},{\bf p}_0)\to
\infty}\frac{\sup_{\xi_C\in
\partial_C^\Delta {\bf f}({\bf p})}\langle \xi_C, \exp_{\bf p}^{-1}({\bf p}_0)\rangle_{\bf g}+
\langle \xi_C^0, \exp_{{\bf p}_0}^{-1}({\bf p})\rangle_{\bf g}}{{\bf
d_g}({\bf p},{\bf p}_0)}< -\|\xi_C^0\|_{\bf g}, \ {\bf p}\in {\bf
K}.$$}

\begin{remark}\rm (a) A similar assumption to hypothesis $(H_{{\bf p}_0})$ can be found in
N\'emeth \cite{Nemeth} in the context of variational inequalities.

(b) Note that for the above numerical example, $(H_{{\bf p}_0})$ is
not satisfied for any ${\bf p}_0=(x_0,y_0)\in [0,\infty)\times
[0,\infty)$. Indeed, one has $L_{(x_0,y_0)}=-e^{x_0+y_0},$ and
$\|\xi_C^0\|_{\bf g}=e^{x_0+y_0}\sqrt{2}$. Therefore, the facts that
${\mathcal S}_{NS}({\bf f,K})={\mathcal S}_{NE}({\bf
f,K})=\emptyset$ are not unexpected.
\end{remark}

The precise statement of the existence result is as follows.

\begin{theorem}\label{begle-noncompact}
Let $(M_i,g_i)$ be finite-dimensional Hadamard manifolds;
$K_i\subset M_i$ be non-empty, closed, geodesic convex sets; and
$U_i\subset M_i$ be open sets containing $K_i$, $i\in \{1,... ,n\}$.
Assume that ${\bf f}\in \mathcal L_{({\bf K,U,M})}$, the map ${\bf
K}\ni {\bf p}\mapsto\partial_C^\Delta {\bf f}({\bf p})$ is upper
semicontinuous, and hypothesis $(H_{{\bf p}_0})$ holds for some
${\bf p}_0\in {\bf K}$.
 Then there exists at
least one Nash-Stampacchia equilibrium point for ${\bf (f,K)}$,
i.e., ${\mathcal S}_{NS}{ \bf (f,K)}\neq \emptyset$.
\end{theorem}
\noindent {\it Proof.} Let $E_0\in \mathbf R$ such that $L_{{\bf
p}_0}<-E_0<-\|\xi_C^0\|_{\bf g}.$ On account of hypothesis $(H_{{\bf
p}_0})$ there exists $R>0$ large enough such that for every ${\bf
p}\in {\bf K}$ with ${\bf d_g}({\bf p},{\bf p}_0)\geq R$, we have
$$\sup_{\xi_C\in
\partial_C^\Delta {\bf f}({\bf p})}\langle \xi_C, \exp_{\bf p}^{-1}({\bf p}_0)\rangle_{\bf g}+
\langle \xi_C^0, \exp_{{\bf p}_0}^{-1}({\bf p})\rangle_{\bf g}\leq
-E_0{\bf d_g}({\bf p},{\bf p}_0).$$ It is clear that ${\bf K}\cap
\overline B_{\bf g}({\bf p}_0,R)\neq \emptyset,$ where $\overline
B_{\bf g}({\bf p}_0,R)$ denotes the closed geodesic ball in ${\bf
(M,g)}$ with center $ {\bf p}_0$ and radius $R.$ In particular, from
(\ref{exp-dist}) and (\ref{cauchy}), for every ${\bf p}\in {\bf K}$
with ${\bf d_g}({\bf p},{\bf p}_0)\geq R$, the above relation yields
\begin{eqnarray}\label{coer-egy-1}
  \sup_{\xi_C\in
\partial_C^\Delta {\bf f}({\bf p})}\langle \xi_C, \exp_{\bf p}^{-1}({\bf p}_0)\rangle_{\bf
g} &\leq& -E_0{\bf d_g}({\bf p},{\bf
p}_0)+\|\xi_C^0\|_{\bf g}\|\exp_{{\bf p}_0}^{-1}({\bf p})\|_{\bf g} \\
 \nonumber  &=& (-E_0+\|\xi_C^0\|_{\bf g}){\bf d_g}({\bf p},{\bf
p}_0)\\\nonumber &<&0.
\end{eqnarray}
Let ${\bf K}_R={\bf K}\cap \overline B_{\bf g}({\bf p}_0,R).$  It is
clear that ${\bf K}_R$ is a geodesic convex, compact subset of ${\bf
M}$. By applying Theorem \ref{begle}, we immediately have that
${\tilde {\bf p}}\in {\mathcal S}_{NS}({\bf f,K}_R)\neq \emptyset$,
i.e., there exists $  \tilde \xi_C\in
\partial_C^\Delta {\bf f}({\tilde {\bf p}})$ such that
\begin{equation}\label{non-comp-egy-2}
\langle \tilde \xi_C,\exp_{{\tilde {\bf p}}}^{-1}({\bf
p})\rangle_{\bf g}\geq 0\ \ {\rm for\ all}\ {\bf p}\in {\bf K}_R.
\end{equation}
It is also clear that ${\bf d_g}({\tilde {\bf p}},{\bf p}_0)< R.$
Indeed, assuming the contrary, we obtain from (\ref{coer-egy-1})
that $\langle \tilde \xi_C,\exp_{{\tilde {\bf p}}}^{-1}({\bf
p}_0)\rangle_{\bf g}< 0$, which contradicts relation
(\ref{non-comp-egy-2}). Now, fix ${\bf q\in K}$ arbitrarily. Thus,
for $\varepsilon > 0$ small enough, the element ${\bf
p}=\exp_{{\tilde {\bf p}}}(\varepsilon \exp_{{\tilde {\bf
p}}}^{-1}({\bf q}))$ belongs both to ${\bf K}$ and $\overline B_{\bf
g}({\bf p}_0,R)$, so ${\bf K}_R.$ By substituting ${\bf p}$ into
(\ref{non-comp-egy-2}), we obtain that $\langle \tilde
\xi_C,\exp_{{\tilde {\bf p}}}^{-1}({\bf q})\rangle_{\bf g}\geq 0$.
The arbitrariness of ${\bf q\in K}$ shows that ${\tilde {\bf p}}\in
{\bf K}$ is actually a Nash-Stampacchia equilibrium point for ${\bf
(f,K)}$, which ends the proof.
 \hfill $\diamondsuit$

\begin{remark}\rm
A relevant application to Theorem \ref{begle-noncompact} will be
given in Example \ref{example-3}.
\end{remark}

The second result in the non-compact case is based on a suitable
Lipschitz-type assumption. In order to avoid technicalities in our
further calculations, we will consider that ${\bf f}\in \mathcal
C_{({\bf K, U,M})}$. In this case,
 $\partial_C^\Delta {\bf f}({\bf p})$ and $A_\alpha^{\bf f}({\bf
p})$ are singletons for every ${\bf p\in K}$ and $\alpha>0.$

For ${\bf f}\in \mathcal C_{({\bf K, U, M})}$, $\alpha>0$ and
$0<\rho<1$  we introduce the hypothesis:\\
 $$(H_{\bf K}^{\alpha,\rho})\ \ {\bf d_g}(\exp_{\bf p}(-\alpha \partial_C^\Delta {\bf f}({\bf p})),\exp_{\bf
 q}(-\alpha
\partial_C^\Delta {\bf f}({\bf q})))\leq (1-\rho) {\bf d_g}({\bf p},{\bf
q}) \ {\rm for\ all}\ {\bf p},{\bf q}\in {\bf K}.$$


\vspace{0.1cm}

\begin{remark}\rm One
can show that $(H_{\bf K}^{\alpha,\rho})$ implies $(H_{{\bf p}_0})$
for every ${\bf p}_0\in {\bf K}$ whenever $(M_i,g_i)$ are Euclidean
spaces. However, it is not clear if the same holds for Hadamard
manifolds.
 \end{remark}

 Finding fixed points for $A_\alpha^{\bf f}$,  one
could expect to apply
{dynamical systems};  we consider both {\it discrete} and {\it continuous} ones. 
First, for some $\alpha>0$ and ${\bf p}_0\in {\bf M}$ fixed, we
consider the discrete dynamical system $$\leqno{(DDS)_\alpha}\ \ \ \
\ {\bf p}_{k+1}= A_\alpha^{\bf f}(P_{\bf K}({\bf p}_k)).$$ Second,
according to Theorem \ref{theorem-equivalence}, we clearly have that
$${\bf p}\in {\mathcal S}_{NS}({\bf f,K}) \Leftrightarrow
0=\exp_{\bf p}^{-1}(A_\alpha^{\bf f}({\bf p}))\ \ {\rm for\
all/some}\ \alpha>0.$$ Consequently, for some $\alpha>0$ and ${\bf
p}_0\in {\bf M}$ fixed, the above equivalence motivates the study of
the continuous dynamical system $$\leqno{(CDS)_\alpha} \ \ \ \ \
\left\{
\begin{array}{lll}
\dot\eta(t)= \exp_{\eta(t)}^{-1}(A_\alpha^{\bf f}(P_{\bf
K}(\eta(t))))
\\ \eta(0)={\bf
p}_0.
\end{array}
\right. $$

\noindent The next result describes the exponential stability of the
orbits in both cases.

\begin{theorem}\label{dyn-system-theorem-fo}
Let $(M_i,g_i)$ be finite-dimensional Hadamard manifolds;
$K_i\subset M_i$ be non-empty, closed geodesics convex sets;
$U_i\subset M_i$ be open sets containing $K_i$; and $f_i:{\bf K}\to
\mathbf R$ be functions, $i\in \{1,... ,n\}$ such that ${\bf f}\in
\mathcal C_{({\bf K,U,M})}.$ Assume that $(H_{\bf K}^{\alpha,\rho})$
holds true for some $\alpha>0$ and $0<\rho<1$. Then the set of
Nash-Stampacchia equilibrium points for ${\bf (f,K)}$ is a
singleton, i.e., ${\mathcal S}_{NS}({\bf f,K}) =\{{\bf \tilde p}\}$.
Moreover, for each ${\bf p}_0\in {\bf M}$, we have
\begin{itemize}
  \item[\rm (i)] the  orbit $\{{\bf p}_k\}$ of $(DDS)_\alpha$ converges
exponentially to ${\bf \tilde p\in K}$ and $${\bf d_g}({\bf
p}_{k},{\bf \tilde p})\leq \frac{(1-\rho)^k}{\rho}{\bf d_g}({\bf
p}_1,{\bf
 p}_0)\ for\ all\ k\in \mathbf N;$$
  \item[\rm (ii)]  the
orbit $\eta$ of $(CDS)_\alpha$ is globally defined on $[0,\infty)$
and it converges exponentially to ${\bf \tilde p\in K}$ and $${\bf
d_g}(\eta(t),{\bf \tilde p})\leq {e}^{-\rho t}{\bf d_g}({\bf
p}_0,{\bf \tilde p})\ for\ all\ t\geq 0.$$
\end{itemize}
Furthermore, the set ${\bf K}$ is invariant with respect to the
orbits in both cases whenever ${\bf p}_0\in {\bf K}$.
\end{theorem}

\noindent {\it Proof.}   Let ${\bf p, q\in M}$ be arbitrarily fixed.
On account of the non-expansiveness of the projection $P_{\bf K}$
(see Proposition \ref{egyik-irany} (ii)) and hypothesis  $(H_{\bf
K}^{\alpha,\rho})$, we have
that\vspace{0.2cm}\\
${\bf d_g}((A_\alpha^{\bf f}\circ P_{\bf K})({\bf p}),(A_\alpha^{\bf
f}\circ P_{\bf K})({\bf q}))$
\begin{eqnarray*}
&=& {\bf d_g}(P_{\bf K}(\exp_{P_{\bf K}({\bf
p})}(-\alpha\partial_{C}^\Delta {\bf f}(P_{\bf K}({\bf p})))),P_{\bf
K}(\exp_{P_{\bf K}({\bf q})}(-\alpha\partial_{C}^\Delta {\bf
f}(P_{\bf K}({\bf q})))))\\&\leq& {\bf d_g}(\exp_{P_{\bf K}({\bf
p})}(-\alpha\partial_{C}^\Delta {\bf f}(P_{\bf K}({\bf
p}))),\exp_{P_{\bf K}({\bf q})}(-\alpha\partial_{C}^\Delta {\bf
f}(P_{\bf K}({\bf q}))))\\&\leq&(1-\rho){\bf d_g}(P_{\bf K }({\bf
p}),P_{\bf K}({\bf q}))\\&\leq&(1-\rho){\bf d_g}({\bf p},{\bf q}),
\end{eqnarray*}
which means that the map $A_\alpha^{\bf f}\circ P_{\bf K}:{\bf M\to
M}$ is a $(1-\rho)$-contraction on ${\bf M}$.

(i) 
Since $({\bf M,d_g})$ is a complete metric space, 
a standard Banach fixed point argument shows that $A_\alpha^{\bf
f}\circ P_{\bf K}$ has a unique fixed point $\tilde {\bf p}\in M.$
Since Im$A_\alpha^{\bf f}\subset {\bf K}$, then ${\bf \tilde p}\in
{\bf K}$. Therefore, we have that $A_\alpha^{\bf f}({\bf \tilde
p})={\bf \tilde p}.$ Due to Theorem \ref{theorem-equivalence},
${\mathcal S}_{NS}({\bf f,K})=\{{\bf \tilde p}\}$
 and the estimate for ${\bf d_g}({\bf
p}_k,{\bf \tilde p})$ yields in a usual manner.

(ii) 
Since $A_\alpha^{\bf f}\circ P_{\bf K}:{\bf M\to M}$ is a
$(1-\rho)$-contraction on ${\bf M}$ (thus locally Lipschitz in
particular),  the map ${\bf M}\ni {\bf p}\mapsto G({\bf
p}):=\exp_{\bf p}^{-1}(A_\alpha^{\bf f}(P_{\bf K}({\bf p})))$ is of
class $C^{1-0}$. Now, we may guarantee the existence of a unique
maximal orbit $\eta:[0,T_{\rm max})\to {\bf M}$ of $(CDS)_\alpha.$

We assume that $T_{\rm max}<\infty.$  Let us consider the Lyapunov
function $h:[0,T_{\rm max})\to \mathbf R$  defined by
$$h(t)=\frac{1}{2}{\bf d}_{\bf g}^2(\eta(t),{\bf \tilde p}).$$ The
function $h$ is differentiable for a.e. $t\in [0,T_{\rm max})$ and
in the differentiable points of $\eta$ we have
\begin{eqnarray*}
h'(t)&=&-{\bf g}(\dot \eta(t),\exp_{\eta(t)}^{-1}({\bf \tilde p}))\\
&=&-{\bf g}(\exp_{\eta(t)}^{-1}(A_\alpha^{\bf f}(P_{\bf
K}(\eta(t)))),\exp_{\eta(t)}^{-1}({\bf \tilde p}))\ \ \ \ \ \ \ \ \
\ \ \ \ \ \ \ \ \ \ \ \ \ \ \ \ \ \ ({\rm cf.}\ {(CDS)_\alpha})
\\ &=& -{\bf g}(\exp_{\eta(t)}^{-1}(A_\alpha^{\bf f}(P_{\bf K}(\eta(t))))-\exp_{\eta(t)}^{-1}({\bf \tilde
p}),\exp_{\eta(t)}^{-1}({\bf \tilde p}))\\&& -{\bf
g}(\exp_{\eta(t)}^{-1}({\bf \tilde p}),\exp_{\eta(t)}^{-1}({\bf
\tilde p}))\\ &\leq & \|\exp_{\eta(t)}^{-1}(A_\alpha^{\bf f}(P_{\bf
K}(\eta(t))))-\exp_{\eta(t)}^{-1}({\bf \tilde p})\|_{\bf
g}\cdot\|\exp_{\eta(t)}^{-1}({\bf \tilde p}))\|_{\bf
g}-\|\exp_{\eta(t)}^{-1}({\bf \tilde p}))\|_{\bf g}^2.
\end{eqnarray*}
In the last estimate we used the Cauchy-Schwartz inequality
(\ref{cauchy}). From (\ref{exp-dist}) we have that
\begin{equation}\label{tavolsag-norm}
  \|\exp_{\eta(t)}^{-1}({\bf \tilde p}))\|_{\bf g}={\bf d_g}(\eta(t),{\bf \tilde
  p}).
\end{equation}
We claim that for every $t\in [0,T_{\rm
max})$ one has {
\begin{equation}\label{nagy-kerdes-Rauch}
\|\exp_{\eta(t)}^{-1}(A_\alpha^{\bf f}(P_{\bf
K}(\eta(t))))-\exp_{\eta(t)}^{-1}({\bf \tilde p})\|_{\bf g}\leq {\bf
d_g}(A_\alpha^{\bf f}(P_{\bf K}(\eta(t))),{\bf \tilde p}).
\end{equation}}
 To see this, fix a point   $t\in [0,T_{\rm max})$ where $\eta$ is differentiable, and  let
$\gamma:[0,1]\to {\bf M}$, $\tilde \gamma:[0,1]\to T_{\eta(t)}{\bf
M}$ and $\overline \gamma:[0,1]\to T_{\eta(t)}{\bf M}$ be three
curves such that
\begin{itemize}
  \item[$\bullet$] $\gamma$ is the unique minimal geodesic joining the two points $\gamma(0)={\bf \tilde p}\in
{\bf K}$ and $\gamma(1)=A_\alpha^{\bf f}(P_{\bf K}(\eta(t)))$;
  \item[$\bullet$] $\tilde \gamma(s)=\exp_{\eta
(t)}^{-1}(\gamma(s)),$ $s\in [0,1];$
  \item[$\bullet$] $\overline\gamma(s)=(1-s)\exp_{\eta
(t)}^{-1}({\bf \tilde p})+s\exp_{\eta (t)}^{-1}(A_\alpha^{\bf
f}(P_{\bf K}(\eta(t)))),$ $s\in [0,1].$
\end{itemize}
By the definition of $\gamma,$ we have that
\begin{equation}\label{hossz-1}
 L_{\bf g} (\gamma)={\bf d_g}(A_\alpha^{\bf f}(P_{\bf K}(\eta(t))),{\bf \tilde p}).
\end{equation}
Moreover, since $\overline \gamma$ is a segment of the straight line
in $T_{\eta(t)}{\bf M}$ that joins the endpoints of $\tilde \gamma$,
we have that
\begin{equation}\label{hossz-2}
  l(\overline \gamma)\leq l(\tilde \gamma).
\end{equation}
Here, $l$ denotes the length function on $T_{\eta(t)}{\bf M}$.
Moreover, since the curvature of $({\bf M,g})$ is non-positive, we
may apply a Rauch-type comparison result for the lengths of $\gamma$
and $\tilde \gamma$, see do Carmo \cite[Proposition 2.5,
p.218]{doCarmo}, obtaining that
\begin{equation}\label{hossz-3}
  l(\tilde \gamma)\leq  L_{\bf g}(\gamma).
\end{equation}
Combining relations (\ref{hossz-1}), (\ref{hossz-2}) and
(\ref{hossz-3}) with the fact that $$l(\overline
\gamma)=\|\exp_{\eta (t)}^{-1}(A_\alpha^{\bf f}(P_{\bf K}(\eta(t))))
-\exp_{\eta(t)}^{-1}({\bf \tilde p})\|_{\bf g},$$ relation
(\ref{nagy-kerdes-Rauch}) holds true.

Coming back to $h'(t)$, in view of (\ref{tavolsag-norm}) and
(\ref{nagy-kerdes-Rauch}), it turns out that
\begin{equation}\label{utolso-ki-tudja}
   h'(t)\leq {\bf d_g}(A_\alpha^{\bf f}(P_{\bf K}(\eta(t))),{\bf \tilde p})\cdot{\bf d_g}(\eta(t),{\bf \tilde
  p})-{\bf d}_{\bf g}^2(\eta(t),{\bf \tilde
  p}).
\end{equation}
 On the other
hand, note that  ${\bf \tilde p}\in {\mathcal S}_{NS}({\bf f,K})$,
i.e., $A_\alpha^{\bf f}({\bf \tilde p})={\bf \tilde p}.$ By
exploiting the non-expansiveness of the projection operator $P_{\bf
K}$, see Proposition \ref{egyik-irany} (ii), and $(H_{\bf
K}^{\alpha,\rho})$, we have that
\begin{eqnarray*}
{\bf d_g}(A_\alpha^{\bf f}(P_{\bf K}(\eta(t))),{\bf \tilde p})&=&{\bf d_g}(A_\alpha^{\bf f}(P_{\bf K}(\eta(t))),A_\alpha^{\bf f}({\bf \tilde p}))\\
&=&{\bf d_g}(P_{\bf K}(\exp_{P_{\bf K}(\eta(t))}(-\alpha\partial
_{C}^\Delta{\bf f }(P_{\bf K}(\eta(t))))),P_{\bf K}(\exp_{{\bf
\tilde p}}(-\alpha\partial _{C}^\Delta{\bf f }({\bf \tilde p}))))\\
&\leq& {\bf d_g}(\exp_{P_{\bf K}(\eta(t))}(-\alpha\partial
_{C}^\Delta{\bf f }(P_{\bf K}(\eta(t)))),\exp_{{\bf \tilde
p}}(-\alpha\partial _{C}^\Delta{\bf f }({\bf \tilde p})))\\&\leq&
(1-\rho){\bf d_g}(P_{\bf K}(\eta(t)), {\bf \tilde
p})\\&=&(1-\rho){\bf d_g}(P_{\bf K}(\eta(t)),P_{\bf K}( {\bf \tilde
p}))\\&\leq& (1-\rho){\bf d_g}(\eta(t),{\bf \tilde
p}).\end{eqnarray*} Combining the above relation with
(\ref{utolso-ki-tudja}), for a.e. $t\in [0,T_{\rm max})$ it yields
\begin{eqnarray*}
 h'(t)&\leq&
(1-\rho){\bf d}_{\bf g}^2(\eta(t),{\bf \tilde p})-{\bf d}_{\bf
g}^2(\eta(t),{\bf \tilde p})=-\rho{\bf d}_{\bf g}^2(\eta(t),{\bf
\tilde p}),
\end{eqnarray*}
 which is nothing but $$ h'(t)\leq -2\rho h(t)\ \ \mbox{for a.e.}\ t\in [0,T_{\rm max}). $$ Due
to the latter inequality, we have that $$\frac{d}{dt}[h(t){
e}^{2\rho t}]=[h'(t)+2\rho h(t)]{e}^{2\rho t}\leq 0\ \ \mbox{for
a.e.}\ t\in [0,T_{\rm max}).$$ After integration, one gets
\begin{equation}\label{csokken-novekszik}
  h(t){e}^{2\rho t}\leq h(0)\ \ \mbox{for all}\ t\in [0,T_{\rm max}).
\end{equation}
According to (\ref{csokken-novekszik}), the function $h$ is bounded
on $[0,T_{\rm max})$; thus, there exists ${\bf \overline p}\in {\bf
M}$ such that
$  \lim_{t\nearrow T_{\rm max}}\eta(t)={\bf \overline p}.$
The last limit means that $\eta$ can be extended toward the value
$T_{\max}$,  which  contradicts the maximality of $T_{\rm max}$.
Thus, $T_{\rm max}=\infty.$

Now, relation (\ref{csokken-novekszik}) leads to the required
estimate; indeed, we have  $${\bf d}_{\bf g}(\eta(t),{\bf \tilde
p})\leq {e}^{-\rho t}{\bf d}_{\bf g}(\eta(0),{\bf \tilde p})={
e}^{-\rho t}{\bf d}_{\bf g}({\bf p}_0,{\bf \tilde p})\ \ \mbox{for
all}\ t\in [0,\infty),$$ which concludes the proof of (ii).

Now, we assume that ${\bf p}_0\in {\bf K}$ and we are dealing with
the viability of the orbits for problems $(DDS)_\alpha$ and
$(CDS)_\alpha$, respectively. First, since Im$A_\alpha^{\bf
f}\subset {\bf K}$, then the orbit of $(DDS)_\alpha$ belongs to
${\bf K}$, i.e., ${\bf p}_k\in {\bf K}$ for every $k\in \mathbf N.$
Second, we shall prove that ${\bf K}$ is invariant with respect to
the solutions of $(CDS)_\alpha$, i.e., the image of the global
solution
$\eta:[0,\infty)\to {\bf M}$ of $(CDS)_\alpha$ with $\eta(0)={\bf p}_0\in {\bf K}$, 
entirely belongs to the set ${\bf K}$. To show the latter fact, we
are going to apply Proposition \ref{weak-inva-mas} by choosing
$M:={\bf M}$ and $G:{\bf M}\to T{\bf M}$ defined by  $G({\bf
p}):=\exp_{\bf p}^{-1}(A_\alpha^{\bf f}(P_{\bf K}({\bf p})))$. 

Fix ${\bf p\in K}$ and $\xi\in N_F({\bf p;K})$. Since ${\bf K}$ is
geodesic convex in $({\bf M,g})$, on account of Theorem
\ref{th-normal-cones-tetel}, we have that $\langle \xi, \exp_{\bf
p}^{-1}({\bf q}) \rangle_{\bf g}\leq 0$ for all ${\bf q\in K}.$
In particular, if we choose ${\bf q}_0=A_\alpha^{\bf f}(P_{\bf K}({\bf p}))\in {\bf K},$ 
it turns out that $$H_G({\bf p},\xi)= \langle\xi, G({\bf
p})\rangle_{\bf g}=\langle\xi, \exp_{\bf p}^{-1}(A_\alpha^{\bf
f}(P_{\bf K}({\bf p})))\rangle_{\bf g}=\langle \xi, \exp_{\bf
p}^{-1}({\bf q}_0) \rangle_{\bf g}\leq 0.$$ Consequently, our claim
is proved by applying Proposition \ref{weak-inva-mas}. \hfill
$\diamondsuit$

\section{Curvature rigidity: metric projections versus Hadamard manifolds}\label{section-Nash-curvature}

The obtuse-angle property and the non-expansi\-ve\-ness of $P_{\bf
K}$ for the closed, geodesic convex set ${\bf K\subset M}$ played
indispensable roles in the proof of Theorems
\ref{theorem-equivalence}-\ref{dyn-system-theorem-fo}, which are
well-known features of Hadamard manifolds (see Proposition
\ref{egyik-irany}). In Section \ref{section-main} the product
manifold ${\bf (M,g)}$ is considered to be a Hadamard one due to the
fact that $(M_i,g_i)$ are Hadamard manifolds themselves for each
$i\in \{1,...,n\}$. We actually have the following characterization
which is also of geometric interests in its own right and entitles
us to assert that Hadamard manifolds are the natural framework to
develop the theory of Nash-Stampacchia equilibria on manifolds.

\begin{theorem}\label{fotetel-1} Let $(M_i,g_i)$ be complete, simply connected
Riemannian manifolds, $i\in \{1,...,n\}$, and ${\bf (M,g)}$ their
product manifold. The following statements are equivalent:
\begin{enumerate}
\item[{\rm (i)}] Any non-empty, closed, geodesic convex set $
{\bf K\subset M}$ verifies the obtuse-angle property and $P_{\bf K}$
is non-expansive;

\item[({\rm ii)}]  $(M_i,g_i)$ are Hadamard manifolds for every $i\in \{1,...,n\}$.
\end{enumerate}
\end{theorem}

\noindent {\it Proof.} (ii)$\Rightarrow$(i). As mentioned before, if
$(M_i,g_i)$ are Hadamard manifolds for every $i\in \{1,...,n\}$,
then ${\bf (M,g)}$ is also a Hadamard manifold, see Ballmann
\cite[Example 4, p.147]{Ballmann} and O'Neill \cite[Lemma 40, p.
209]{ONeill}. It remains to apply Proposition \ref{egyik-irany} for
the Hadamard manifold ${\bf (M,g)}$.

(i)$\Rightarrow$(ii). We first prove that ${\bf (M,g)}$ is a
Hadamard manifold. Since $(M_i,g_i)$ are complete and simply
connected Riemannian manifolds for every $i\in \{1,...,n\}$, the
same is true for ${\bf (M,g)}$. We now show that the sectional
curvature of ${\bf (M,g)}$ is non-positive. To see this, let ${\bf
p\in M}$ and ${\bf W}_0,{\bf V}_0\in  T{\bf_pM\setminus \{0\}}.$ We
claim that the sectional curvature of the two-dimensional subspace
$S=$span$\{ {\bf W}_0,{\bf V}_0 \}\subset  T{\bf _pM}$ at the point
${\bf p}$ is non-positive, i.e., $K{\bf _p}(S)\leq 0$. We assume
without loosing the generality that ${\bf V_0}$ and ${\bf W_0}$ are
${\bf g}$-perpendicular, i.e., ${\bf g}({\bf W}_0,{\bf V}_0)=0$.

Let us fix $r_{\bf p}>0$ and $\delta>0$ such that $B{\bf _g}({\bf
p},r_{\bf p})$ is a totally normal ball of ${\bf p}$ and
\begin{equation}\label{totally}
  \delta\left(\|{\bf W}_0\|_{\bf g}+2\|{\bf V}_0\|_{\bf g}\right)<r_{\bf p}.
\end{equation}
Let $\sigma:[-\delta,2\delta]\to {\bf M}$ be the geodesic segment
$\sigma(t)=\exp_{\bf p}(t{\bf V}_0)$ and ${\bf W}$ be the unique
parallel vector field along $\sigma$ with the initial data ${\bf
W}(0)={\bf W}_0$. For any $t\in [0,\delta]$, let
$\gamma_t:[0,\delta]\to {\bf M}$ be the geodesic segment
$\gamma_t(u)=\exp_{\sigma(t)}(u{\bf W}(t)).$

Let us fix $t,u\in [0,\delta]$ arbitrarily, $u\neq 0.$ Due to
(\ref{totally}), the geodesic segment $\gamma_t|_{[0,u]}$ belongs to
the totally normal ball $B{\bf _g}({\bf p},r_{\bf p})$ of ${\bf p}$;
thus, $\gamma_t|_{[0,u]}$ is the unique minimal geodesic joining the
point  $\gamma_t(0)=\sigma(t)$ to $\gamma_t(u).$ Moreover, since
${\bf W}$ is the parallel transport of ${\bf W}(0)={\bf W}_0$ along
$\sigma$, we have ${\bf g}({\bf W}(t),\dot \sigma(t))={\bf g(
W}(0),\dot \sigma(0)) ={\bf g}({\bf W}_0,{\bf V}_0)=0;$ therefore,
$$
  {\bf g}(\dot \gamma_t(0),\dot \sigma(t)) ={\bf g( W}(t),\dot \sigma(t))=0.
$$
Consequently,  the minimal geodesic segment $\gamma_t|_{[0,u]}$
joining $\gamma_t(0)=\sigma(t)$ to $\gamma_t(u)$, and the set ${\bf
K}={\rm Im}\sigma=\{\sigma(t):t\in [-\delta,2\delta]\}$ fulfill
hypothesis $(OA_2)$. Note that ${\rm Im}\sigma$ is a closed,
geodesic convex set in ${\bf M}$; thus, from hypothesis (i) it
follows that ${\rm Im}\sigma$ verifies the obtuse-angle property and
 $P_{{\rm Im}\sigma}$ is non-expansive. Thus, $(OA_2)$
implies $(OA_1)$, i.e., for every $t,u\in [0,\delta]$, we have
$\sigma(t)\in P_{{\rm Im}\sigma}(\gamma_t(u)).$ Since ${\rm
Im}\sigma$ is a Chebyshev set (cf. the non-expansiveness of $P_{{\rm
Im}\sigma}$), for every $t,u\in [0,\delta]$, we have
\begin{equation}\label{proj-egy-csomo}
  P_{{\rm
Im}\sigma}(\gamma_t(u))=\{\sigma(t)\}.
\end{equation}
Thus, for every $t,u\in [0,\delta]$, relation (\ref{proj-egy-csomo})
and the non-expansiveness of $P_{{\rm Im}\sigma}$ imply
\begin{eqnarray}\label{mas-11}
{\bf d_g(p},\sigma(t))& =& {\bf d_g}(\sigma(0),\sigma(t))={\bf d_g}(
P_{{\rm Im}\sigma}(\gamma_0(u)), P_{{\rm Im}\sigma}(\gamma_t(u)))
\\& \leq & {\bf d_g}( \gamma_0(u), \gamma_t(u)).\nonumber
\end{eqnarray}
The above construction (i.e., the parallel transport of ${\bf
W}(0)={\bf W}_0$ along $\sigma$) and the formula of the sectional
curvature in the parallelogramoid of Levi-Civita defined by the
points $p$, $\sigma(t)$, $\gamma_0(u),$ $\gamma_t(u)$ give $$
  K{\bf _p}(S)=\lim_{u,t\to 0}\frac{{\bf d^{\rm 2}_g}({\bf p},\sigma(t))-{\bf d^{\rm 2}_g}(
\gamma_0(u), \gamma_t(u))}{{\bf d_g}({\bf p}, \gamma_0(u))\cdot {\bf
d_g}({\bf p}, \sigma(t))}. $$
 According to
(\ref{mas-11}), the latter limit is non-positive, so
$K{\bf_p}(S)\leq 0,$ which concludes the first part, namely, $({\bf
M,g})$ is a Hadamard manifold.

Now, the main result of Chen \cite[Theorem 1]{CH-Chen-TAMS} implies
that the metric spaces $(M_i,d_{g_i})$ are Aleksandrov NPC spaces
for every $i\in \{1,...,n\}$. Consequently, for each $i\in
\{1,...,n\}$, the Riemannian manifolds $(M_i,g_i)$ have non-positive
sectional curvature, thus they are Hadamard manifolds. The proof is
complete. \hfill $\diamondsuit$

\begin{remark}\rm  The obtuse-angle property and the non-expansiveness of the metric
projection are also key tools behind the theory of monotone vector
fields, proximal point algorithms and variational inequalities
developed on Hadamard manifolds; see Li, L\'opez and Mart\'\i
n-M\'arquez \cite{LLMM-1}, \cite{LLMM-2}, and N\'emeth
\cite{Nemeth}. Within the class of Riemannian manifolds, Theorem
\ref{fotetel-1} shows in particular that Hadamard manifolds are
indeed the appropriate frameworks for developing successfully the
approaches in \cite{LLMM-1}, \cite{LLMM-2}, and \cite{Nemeth} and
further related works.
\end{remark}

\section{Examples}\label{utol-sect}

In this section we present various examples where our main results
can be efficiently applied; for convenience, we give all the details
in our calculations by keeping also the notations from the previous
sections.

\begin{example}\rm \label{example-1}
Let $$K_1=\{(x_1,x_2)\in \mathbf R^2_+: 
x_1^2+x_2^2\leq 4\leq (x_1-1)^2+x_2^2\},\ K_2=[-1,1],$$
and the functions $f_1,f_2:K_1\times K_2 \to \mathbf R$  defined for
$(x_1,x_2)\in K_1$ and $y\in K_2$ by
$$f_1((x_1,x_2),y)=y(x_1^3+y(1-x_2)^3);\ \ f_2((x_1,x_2),y)=-y^2x_2+4|y|(x_1+1).$$
It is clear that $K_1\subset \mathbf R^2$ is not convex in the usual
sense while $K_2\subset \mathbf R$ is. However, if we consider the
Poincar\'e upper-plane model $(\mathbf H^2,g_{\mathbf H})$,  the set
$K_1\subset \mathbf H^2$ is geodesic convex (and compact) with
respect to the metric $g_{\mathbf H}=(\frac{\delta_{ij}}{x_2^2})$.
Therefore, we embed the set $K_1$ into the Hadamard manifold
$(\mathbf H^2,g_{\mathbf H})$, and $K_2$ into the standard Euclidean
space $(\mathbf R,g_0)$. After natural extensions of $f_1(\cdot,y)$
and $f_2((x_1,x_2),\cdot)$ to the whole $U_1=\mathbf H^2$ and
$U_2=\mathbf R$, respectively, we clearly have that $f_1(\cdot,y)$
is a $C^1$ function on $\mathbf H^2$ for every $y\in K_2$, while
$f_2((x_1,x_2),\cdot)$ is a locally Lipschitz function on $\mathbf
R$ for every $(x_1,x_2)\in K_1$. Thus,  ${\bf f}=(f_1,f_2)\in
\mathcal L_{(K_1\times K_2,\mathbf H^2\times \mathbf R,\mathbf
H^2\times \mathbf R)}$ and for every $((x_1,x_2),y)\in {\bf
K}=K_1\times K_2$, we have $$\partial_C^1 f_1((x_1,x_2),y)={\rm
grad}f_1(\cdot,y)(x_1,x_2)=\left(g_{\mathbf H}^{ij}\frac{\partial
f_1(\cdot,y)}{\partial x_j}\right)_i=3yx_2^2(x_1^2,-y(1-x_2)^2);$$
$$\partial_C^2 f_2((x_1,x_2),y)= \left\{
\begin{array}{lll}
-2yx_2-4(x_1+1)& {\rm if} & y<0,\\
4(x_1+1)[-1,1]& {\rm if} & y=0,\\
-2yx_2+4(x_1+1)&  {\rm if} & y>0.
\end{array}
\right. $$ It is now clear that  the map ${\bf K}\ni
((x_1,x_2),y)\mapsto\partial_C^\Delta {\bf f}(((x_1,x_2),y))$ is
upper semicontinuous. Consequently, on account of Theorem
\ref{begle}, ${\mathcal S}_{NS}({\bf f,K})\neq \emptyset,$ and its
elements are precisely the solutions  $((\tilde x_1,\tilde
x_2),\tilde y)\in {\bf K}$ of the  system
$$\left\{
\begin{array}{lll}
\langle\partial_C^1 f_1((\tilde x_1,\tilde x_2),\tilde y),\exp_{(\tilde x_1,\tilde x_2)}^{-1}(q_1,q_2)\rangle_{g_{\mathbf H}}\geq 0 & {\rm \mbox{for all}} &  (q_1,q_2)\in K_1, \\
\xi_C^2(q-\tilde y)\geq 0\ \ {\rm for\ some}\ \xi_C^2\in
\partial_C^2 f_2((\tilde x_1,\tilde x_2),\tilde y)& {\rm \mbox{for
all}} & q\in K_2.
\end{array} \right.\eqno{(S_1)}$$
 In order to solve $(S_1)$ we first observe that
\begin{equation}\label{geom-pelda}
K_1\subset \{(x_1,x_2)\in \mathbf R^2:\sqrt{3}\leq x_2\leq
2(x_1+1)\}.
\end{equation}
 We distinguish four cases:

(a) If $\tilde y=0$ then both inequalities of $(S_1)$ hold for every
$(\tilde x_1,\tilde x_2)\in {K_1}$  by choosing $\xi_C^2=0\in
\partial_C^2 f_2((\tilde x_1,\tilde x_2),0)$ in the second relation.
Thus, $((\tilde x_1,\tilde x_2),0)\in {\mathcal S}_{NS}({\bf f,K})$
for every $(\tilde x_1,\tilde x_2)\in {\bf K}$.

(b) Let $0<\tilde y<1$. The second inequality of $(S_1)$ gives that
$-2\tilde y\tilde x_2+4(\tilde x_1+1)=0$; together with
(\ref{geom-pelda}) it yields $0=\tilde y\tilde x_2-2(\tilde
x_1+1)<\tilde x_2-2(\tilde x_1+1)\leq 0,$ a contradiction.

(c) Let $\tilde y=1$. The second inequality of $(S_1)$ is true if
and only if $-2\tilde x_2+4(\tilde x_1+1)\leq 0$. Due to
(\ref{geom-pelda}), we necessarily have $\tilde x_2=2(\tilde
x_1+1)$; this Euclidean line intersects the set $K_1$ in the unique
point $(\tilde x_1,\tilde x_2)=(0,2)\in K_1.$ By the geometrical
meaning of the exponential map one can conclude that
$$\{t\exp_{(0,2)}^{-1}(q_1,q_2):(q_1,q_2)\in
K_1, t\geq 0\}=\{(x,-y)\in \mathbf R^2:x,y\geq 0\}.$$ Taking into
account this relation and $\partial_C^1 f_1((0,2),1)=(0,-12)$, the
first inequality of $(S_1)$ holds true as well. Therefore,
$((0,2),1)\in {\mathcal S}_{NS}({\bf f,K}).$

(d) Similar reason as in (b) (for $-1<\tilde y<0$) and (c) (for
$\tilde y=-1$) gives that $((0,2),-1)\in {\mathcal S}_{NS}({\bf
f,K}).$
Thus, from (a)-(d) we have that $${\mathcal S}_{NS}({\bf f,K})=
(K_1\times \{0\})\cup \{((0,2),1),((0,2),-1)\}.$$ Now, on account of
Theorem \ref{equi-kritikus} (i) we may easily select the Nash
equilibrium points for ${\bf (f,K)}$ among the elements of
${\mathcal S}_{NS}({\bf f,K})$ obtaining that\\
\hspace*{5cm} ${\mathcal S}_{NE}({\bf f,K})= K_1\times \{0\}.$
\hfill $\diamondsuit$
\end{example}

\vspace{0.4cm}

In the rest of the paper we deal with some applications involving
matrices; thus, we recall some basic notions from the
matrix-calculus.  Fix $n\geq 2$. Let $M_n(\mathbf R)$ be the set of
symmetric $n\times n$ matrices with real values, and $M_n^+(\mathbf
R)\subset M_n(\mathbf R)$ be the cone of symmetric positive definite
matrices. The standard inner product on $M_n(\mathbf R)$ is defined
as
\begin{equation}\label{eukl-metrik-matrix}
    \langle U,V\rangle= {\rm tr}(UV).
\end{equation}
Here, ${\rm tr}(Y)$ denotes the trace of $Y\in M_n(\mathbf R)$. It
is well-known that $(M_n(\mathbf R),\langle \cdot, \cdot \rangle)$
is an Euclidean space, the unique geodesic between $X,Y\in
M_n(\mathbf R)$ is
\begin{equation}\label{geodesic-matrix-eukl}
\gamma_{X,Y}^E(s)=(1-s)X+sY,\ \ s\in [0,1].
\end{equation}

 The set $M_n^+(\mathbf R)$ will be endowed with the Killing form
\begin{equation}\label{hadam-metrika}
  \langle\langle U,V \rangle\rangle_X={\rm tr}(X^{-1}VX^{-1}U),\ \ \ X\in
M_n^+(\mathbf R),\ U,V\in T_X(M_n^+(\mathbf R)).
\end{equation}
Note that the pair $(M_n^+(\mathbf R),\langle\langle\cdot,\cdot
\rangle\rangle)$ is a Hadamard manifold, see Lang \cite[Chapter
XII]{Lang}, and $T_X(M_n^+(\mathbf R))\simeq M_n(\mathbf R).$ The
unique geodesic segment connecting $X,Y\in M_n^+(\mathbf R)$ is
defined by
\begin{equation}\label{geodesic-matrix}
\gamma_{X,Y}^H(s)=X^{1/2}(X^{-1/2}Y X^{-1/2})^s X^{1/2},\ \ s\in
[0,1].
\end{equation}
In particular,
$\frac{d}{ds}{\gamma_{X,Y}^H}(s)|_{s=0}=X^{1/2}\ln(X^{-1/2}Y
X^{-1/2}) X^{1/2};$ consequently, for each $X,Y\in M_n^+(\mathbf
R)$, we have
$$\exp_X^{-1}Y=X^{1/2}\ln(X^{-1/2}Y X^{-1/2}) X^{1/2}.$$ Moreover, the  metric function on $M_n^+(\mathbf R)$ is given by
\begin{equation}\label{dh-metrika}
    d_H^2(X,Y)=\langle\langle
\exp_X^{-1}Y,\exp_X^{-1}Y\rangle\rangle_X={\rm tr}(\ln^2(X^{-1/2}Y
X^{-1/2})).
\end{equation}

\begin{example}\label{example-2}\rm
Let
$$K_1=[0,2],\ K_2=\{X\in M_n^+(\mathbf R):{\rm tr}(\ln^2X)\leq 1\leq \det X \leq 2\},$$
and the functions $f_1,f_2:K_1\times K_2 \to \mathbf R$  defined by
\begin{equation}
\label{f1f2-matrix-1}
  f_1(t,X)=(\max(t,1))^{n-1}{\rm tr}^2(X)-4n\ln (t+1)S_2(X),
\end{equation}
\begin{equation}
\label{f1f2-matrix-2}
  f_2(t,X)=g(t)\left({\rm tr}(X^{-1})+1\right)^{t+1}+h(t) \ln  \det
X.
\end{equation}
Here, $S_2(Y)$ denotes the second elementary symmetric function of
the eigenvalues $\lambda_1,...,\lambda_n$ of $Y$, i.e.,
\begin{equation}\label{s2-szimm}
S_2(Y)=\sum_{1\leq i_1<i_2\leq n}\lambda_{i_1}\lambda_{i_2},
\end{equation}
and $g,h:K_1\to \mathbf R$ are two continuous functions such that
\begin{equation}\label{auxiliar-egyenlte}
    h(t)\geq 2(n+1)g(t)\geq 0\ {\rm for\ all}\ t\in K_1.
\end{equation}

The elements of $\mathcal S_{NE}({\bf f,K})$ are the solutions
$(\tilde t,\tilde X)\in {\bf K}$ of the  system {\small $$\left\{
\begin{array}{lll}
\left[(\max(t,1))^{n-1}-(\max(\tilde t,1))^{n-1}\right]{\rm tr}^2(\tilde X)\geq 4nS_2(\tilde X)\ln \frac{t+1}{\tilde t+1},    & {{\forall}}   t\in K_1, \\
g(\tilde t)\left[({\rm tr}(Y^{-1})+1)^{\tilde t+1}-({\rm tr}(\tilde
X^{-1})+1)^{\tilde t+1}\right]+h(\tilde t) \ln  \frac{\det Y}{\det
\tilde X}\geq 0, & {{\forall}}  Y\in K_2.
\end{array} \right.\eqno{(S_2)}$$}
The involved forms in $(S_2)$ suggest an approach via the
Nash-Stampacchia equilibria for $({\bf f,K});$ first of all, we have
to find the appropriate context where the machinery described in \S
\ref{section-main} works efficiently.

At first glance, the natural geometric framework seems to be
$M_n(\mathbf R)$ with the inner product $\langle\cdot, \cdot\rangle$
defined in (\ref{eukl-metrik-matrix}). Note however that the set
$K_2$ is not geodesic convex with respect to $\langle \cdot, \cdot
\rangle$. Indeed, let $X={\rm diag}(2,1,...,1)\in K_2$ and $Y={\rm
diag}(1,2,...,1)\in K_2$ and $\gamma_{X,Y}^E$ be the Euclidean
geodesic connecting them, see (\ref{geodesic-matrix-eukl}); although
$\gamma_{X,Y}^E(s)\in M_n^+(\mathbf R)$ and ${\rm
tr}(\ln^2(\gamma_{X,Y}^E(s)))=\ln^2(2-s)+\ln^2(1+s)\leq \ln^2 2$ for
every $s\in [0,1]$, we have that $\det(\gamma_{X,Y}^E(s))>2$ for
every $0<s<1.$ Consequently, a more appropriate metric is needed to
provide some sort of geodesic convexity for $K_2.$ To complete this
fact, we restrict our attention to the cone of symmetric positive
definite matrices $M_n^+(\mathbf R)$ with the metric introduced in
(\ref{hadam-metrika}).

Let $I_n\in M_n^+(\mathbf R)$ be the identity matrix, and $\overline
B_H(I_n,1)$ be the closed geodesic ball in $M_n^+(\mathbf R)$ with
center $I_n$ and radius 1. Note that $$K_2=\overline B_H(I_n,1)\cap
\{X\in M_n^+(\mathbf R): 1\leq \det X\leq 2\}.$$ Indeed, for every
$X\in M_n^+(\mathbf R)$, we have
\begin{equation}\label{dist-egyseg-x}
    d_H^2(I_n,X)={\rm tr}(\ln^2 X).
\end{equation}
Since $K_2$ is bounded and closed, on account of the Hopf-Rinow
theorem, $K_2$ is compact. Moreover, as a geodesic ball in the
Hadamard manifold $(M_n^+(\mathbf R),\langle\langle\cdot,\cdot
\rangle\rangle)$, the set $\overline B_H(I_n,1)$ is geodesic convex.
Keeping the notation from (\ref{geodesic-matrix}), if $X,Y\in K_2$,
one has for every $s\in [0,1]$ that
$$\det (\gamma_{X,Y}^H(s))=(\det X)^{1-s}(\det Y)^s\in [1,2],$$
which shows the geodesic convexity of $K_2$ in $(M_n^+(\mathbf
R),\langle\langle\cdot,\cdot \rangle\rangle)$.

After naturally extending the functions $f_1(\cdot,X)$ and
$f_2(t,\cdot)$ to $U_1=(-\frac{1}{2},\infty)$ and $U_2=M_n^+(\mathbf
R)$ by the same expressions (see (\ref{f1f2-matrix-1}) and
(\ref{f1f2-matrix-2})), we clearly have that ${\bf f}=(f_1,f_2)\in
\mathcal L_{({\bf K, U, M})}$, where ${\bf U}=U_1\times U_2,$ and
${\bf M}=\mathbf R\times M_n^+(\mathbf R)$. A standard computation
shows that for every $(t,X)\in U_1\times K_2$, we have
$$\partial_C^1 f_1( t, X)
=-\frac{4nS_2( X)}{ t+1}+{\rm tr}^2( X)\cdot\left\{
\begin{array}{lll}
0 & {\rm \mbox{if}} &  -1/2<  t<1, \\
{[0,n-1]} & {\rm \mbox{if}} &  t =1,\\
(n-1) t^{n-2} & {\rm \mbox{if}} & 1< t.
\end{array} \right.$$
For every $t\in K_1$, the Euclidean gradient of $f_2(t,\cdot)$ at
$X\in U_2=M_n^+(\mathbf R)$ is $$f_2'( t,\cdot)(
X)=-g(t)(t+1)\left({\rm tr}(X^{-1})+1\right)^{t}X^{-2}+h(t)
 X^{-1},$$ thus  the Riemannian derivative has the form
\begin{eqnarray*}
  \partial_C^2 f_2(t, X) &=& {\rm grad} f_2( t,\cdot)( X)= X
f_2'( t,\cdot)( X) X \\
   &=&-g(t)(t+1)\left({\rm
tr}(X^{-1})+1\right)^{t}I_n+h(t)
 X.
\end{eqnarray*}
The above expressions show that ${\bf K}\ni
(t,X)\mapsto\partial_C^\Delta {\bf f}(t,X)$ is upper semicontinuous.
Therefore, Theorem \ref{begle} implies that ${\mathcal S}_{NS}({\bf
f,K})\neq \emptyset,$ and its elements $(\tilde t,\tilde X)\in {\bf
K}$ are precisely the solutions of the system
$$\left\{
\begin{array}{lll}
\xi_1(t-\tilde t)\geq 0 \ \ {\rm for\ some}\ \ \xi_1\in \partial_C^1 f_1(\tilde t,\tilde X) & {\rm \mbox{for all}} &  t\in K_1, \\
 \langle\langle \partial_C^2 f_2(\tilde t,\tilde X),\exp^{-1}_{\tilde
X}Y\rangle\rangle_{\tilde X}\geq 0& {\rm \mbox{for all}} & Y\in K_2,
\end{array} \right.\eqno{(S_2')}$$
We notice that the solutions of $(S_2')$ and $(S_2)$ coincide. In
fact, we may show that ${\bf f}\in \mathcal K_{({\bf K, U, M})}$;
thus from Theorem \ref{equi-kritikus} (ii) we have that $\mathcal
S_{NE}({\bf f,K})=\mathcal S_{NS}({\bf f,K})=\mathcal S_{NC}({\bf
f,K})$. It is clear that the map $t\mapsto f_1(t,X)$ is convex on
$U_1$ for every $X\in K_2$. Moreover, $X\mapsto f_2(t,X)$ is also a
convex function on $U_2=M_n^+(\mathbf R)$ for every $t\in K_1$.
Indeed, fix
 $X,Y\in K_2$ and let $\gamma_{X,Y}^H:[0,1]\to
K_2$ be the unique geodesic segment connecting $X$ and $Y$, see
(\ref{geodesic-matrix}). For every $s\in [0,1]$, we have that
\begin{eqnarray*}
\ln \det  (\gamma_{X,Y}^H(s))&=& \ln ((\det X)^{1-s}(\det Y)^{s})
\\& = &
(1-s)\ln \det X+s\ln\det Y
\\& = & (1-s)\ln \det (\gamma_{X,Y}^H(0))+s\ln\det (\gamma_{X,Y}^H(1)).\nonumber
\end{eqnarray*}
The Riemannian Hessian of $X\mapsto {\rm tr}(X^{-1})$ with respect
to $\langle\langle\cdot, \cdot \rangle\rangle$ is $${\rm Hess}({\rm
tr}(X^{-1}))(V,V)={\rm tr}(X^{-2}V X^{-1}V
)=|X^{-1}VX^{-1/2}|_F^2\geq 0,
$$
where $|\cdot|_F$ denotes the standard Fr\"obenius norm. Thus,
$X\mapsto {\rm tr}(X^{-1})$ is convex (see Udri\c ste \cite[\S
3.6]{Udriste}), so $X\mapsto({\rm tr}(X^{-1})+1)^{t+1}.$ Combining
the above facts with the non-negativity of $g$ and $h$ (see
(\ref{auxiliar-egyenlte})), it yields that ${\bf f}\in \mathcal
K_{({\bf K, U, M})}$ as we claimed.

By recalling the notation from (\ref{s2-szimm}), the inequality of
Newton has the form
\begin{equation}\label{newton}
S_2(Y)\leq \frac{n-1}{2n}{\rm tr}^2(Y)\ \ {\rm \mbox{for all}}\ \
Y\in M_n(\mathbf R).
\end{equation}
The  possible cases are as follow:

(a) Let $0\leq \tilde t<1$. Then the first relation from $(S_2')$
implies $-\frac{4nS_2(\tilde X)}{\tilde t+1}\geq 0,$ a
contradiction.

(b) If $1<\tilde t<2$,  the first inequality from $(S_2')$ holds if
and only if  $$S_2(\tilde X)=\frac{n-1}{4n}\tilde t^{n-2}(\tilde
t+1){\rm tr}^2(\tilde X),$$ which contradicts Newton's inequality
(\ref{newton}).

(c) If $\tilde t=2$,  from the first inequality of $(S_2')$ it
follows that $$3(n-1)2^{n-4}{\rm tr}^2(\tilde X)\leq {nS_2(\tilde
X)},$$ contradicting again (\ref{newton}).

(d ) Let $\tilde t=1$. From the first relation of $(S_2')$ we
necessarily have that $0=\xi_1\in \partial_C^1 f_1(1,\tilde X).$
This fact is equivalent to $$\frac{2n S_2(\tilde X)}{{\rm
tr}^2(\tilde X)}\in [0,n-1],$$ which holds true, see (\ref{newton}).
In this case, the second relation from $(S_2')$ becomes
$$
-2g(1)\left({\rm tr}(\tilde X^{-1})+1\right)\langle\langle
I_n,\exp^{-1}_{\tilde X}Y\rangle\rangle_{\tilde X}+h(1)
 \langle\langle\tilde X,\exp^{-1}_{\tilde
X}Y\rangle\rangle_{\tilde X}\geq 0,\ \forall Y\in K_2.
$$
By using (\ref{hadam-metrika}) and the well-known formula $e^{{\rm
tr}(\ln X)}= \det X$, the above inequality reduces to {\small
\begin{equation}\label{s2-ujra-masodik} -2g(1)({\rm tr}(\tilde
X^{-1})+1){\rm tr}(\tilde X^{-1}\ln (\tilde X^{-1/2}Y\tilde
X^{-1/2}))+h(1)
 \ln \frac{\det Y}{\det \tilde X}\geq 0,\ \forall Y\in K_2.
\end{equation}}
We also distinguish three cases:

(d1) If $g(1)=h(1)=0$, then $\mathcal S_{NE}({\bf f,K})=\mathcal
S_{NS}({\bf f,K})=\{1\}\times K_2.$

(d2) If $g(1)=0$ and $h(1)>0$, then (\ref{s2-ujra-masodik}) implies
that $\mathcal S_{NE}({\bf f,K})=\mathcal S_{NS}({\bf
f,K})=\{(1,\tilde X)\in {\bf K}: \det \tilde X=1\}.$

(d3) If $g(1)>0$, then (\ref{auxiliar-egyenlte}) implies that
$(1,I_n)\in \mathcal S_{NE}({\bf f,K})=\mathcal S_{NS}({\bf f,K}).$
\hfill $\diamondsuit$
\end{example}

\begin{remark}\rm
We easily observed in the case (d3) that $\tilde X=I_n$ solves
(\ref{s2-ujra-masodik}). Note that the same is not evident at all
for the second inequality in $(S_2)$. We also notice that the
determination of the whole set $\mathcal S_{NS}({\bf f,K})$ in (d3)
is quite difficult; indeed, after a simple matrix-calculus we
realize that (\ref{s2-ujra-masodik}) is equivalent to the equation
$$\tilde X=P_{K_2}\left(e^{-\frac{h(1)}{2g(1)({\rm tr}(\tilde
X^{-1})+1)}}\tilde X e^{\tilde X^{-1}}\right),$$ where $P_{K_2}$ is
the metric projection with respect to the metric $d_H$.
\end{remark}

\begin{example}\rm \label{example-3} Let
$$K_1=[0,\infty),\ K_2=\{X\in M_n^+(\mathbf R):{\rm tr}(X^{-1})\leq n\},$$
and the functions $f_1,f_2:K_1\times K_2 \to \mathbf R$  defined by
$$  f_1(t,X)=t^3\det X-(t-1){\rm tr}(X^{-1}),\ f_2(t,X)=g(t){\rm tr}(\ln^2 X)+h(t) {\rm tr}(X^{-1}),$$
where $g,h:K_1\to \mathbf R$ are two continuous functions such that
\begin{equation}\label{auxiliar-egyenlte-2}
\inf_{K_1} g>0\ {\rm and}\ h\ {\rm is\ bounded}.
\end{equation}
Here,  $M_n^+(\mathbf R)$ is endowed with the Riemannian metric
$\langle\langle\cdot, \cdot\rangle\rangle$ defined in
(\ref{hadam-metrika}). Note that $K_2$ is geodesic convex in
$(M_n^+(\mathbf R),\langle\langle\cdot, \cdot\rangle\rangle)$ being
a sub-level set of the convex function $X\mapsto {\rm tr}(X^{-1})$,
see Example \ref{example-2}. However, $K_2$ is not compact in
$(M_n^+(\mathbf R),\langle\langle\cdot, \cdot\rangle\rangle)$.
Indeed, it is clear that $X_k=kI_n\in K_2$ ($k\geq 1$), but on
account of (\ref{dist-egyseg-x}), we have $d_H(I_n,X_k)=\sqrt{n}\ln
k\to \infty$ as $k\to \infty.$

By extending the functions $f_1(\cdot,X)$ to $\mathbf R$ and
$f_2(t,\cdot)$ to $M_n^+(\mathbf R)$ by the same expressions, it
becomes clear that ${\bf f}=(f_1,f_2)\in \mathcal C_{({\bf
K,U,M})},$ where ${\bf U}={\bf M}=\mathbf R\times M_n^+(\mathbf R)$
is the standard product manifold with metric ${\bf g},$ see
(\ref{prod-metric}). Although the functions $X\mapsto {\rm tr}(\ln^2
X)$ and $X\mapsto {\rm tr}(X^{-1})$ are convex, ${\bf f}=(f_1,f_2)$
does not necessarily belong to $ \mathcal K_{({\bf K,U,M})}$ since
$h$ can be sign-changing. A simple calculation based on
(\ref{dist-egyseg-x}) and (\ref{deriv-tavolsag-hadamard}) gives that
$$\partial_C^\Delta {\bf f}(t,X)=(3t^2\det X-{\rm tr}(X^{-1}),-2g(t)\exp_{X}^{-1}I_n-h(t)I_n),$$ which is a continuous map on ${\bf
K}$.

 In order to
apply Theorem \ref{begle-noncompact}, we will verify hypothesis
$(H_{{\bf p}_0})$ with ${\bf p}_0=(0,I_n)\in {\bf K}.$ Fix ${\bf
p}=(t,X)\in {\bf K}$ arbitrarily. Then, we have
$$\langle
\partial_C^\Delta {\bf
f}({\bf p}_0),\exp_{{\bf p}_0}^{-1}({\bf p})\rangle_{\bf
g}=-nt-h(0)\langle\langle
I_n,\exp_{I_n}^{-1}X\rangle\rangle_{I_n}=-nt-h(0)\ln \det X,$$ while
by (\ref{exp-dist}) and (\ref{dist-egyseg-x}), we also obtain
\begin{eqnarray*}
  \langle \partial_C^\Delta {\bf
f}({\bf p}),\exp_{{\bf p}}^{-1}({\bf p}_0)\rangle_{\bf g} &=&
-t(3t^2\det X-{\rm tr}(X^{-1})) \\
   && + \langle\langle
-2g(t)\exp_{X}^{-1}I_n-h(t)I_n,\exp_{X}^{-1}I_n\rangle\rangle_{X}\\
   &=&-3t^3\det X+t{\rm tr}(X^{-1})\\ && -2g(t){\rm tr}(\ln^2
   X)-h(t)\langle\langle I_n,\exp_{X}^{-1}I_n\rangle\rangle_{X}.
\end{eqnarray*}
On the one hand, by (\ref{cauchy}), (\ref{exp-dist}), and the
inequality ${\rm tr}(A^2)\leq {\rm tr}^2(A)$ (for $A\in
M_n^+(\mathbf R))$, we have
\begin{eqnarray*}
  |\langle\langle I_n,\exp_{X}^{-1}I_n\rangle\rangle_{X}| &\leq & \|I_n\|_X\|\exp_{X}^{-1}I_n\|_X=({\rm tr}(X^{-2}))^{1/2}d_H(I_n,X) \\
   &\leq& {\rm tr}(X^{-1})d_H(I_n,X)\leq nd_H(I_n,X)=n({{\rm tr}(\ln^2
   X)})^{1/2}.
\end{eqnarray*}
On the other hand, since $X\in K_2$, one has $$({\det
X^{-1}})^{1/n}\leq \frac{{\rm tr}(X^{-1})}{n}\leq 1,$$ thus, $\det
X\geq 1.$ On account of (\ref{auxiliar-egyenlte-2}) and the above
estimations, one has {\small $$\frac{-3t^3\det X-2g(t){\rm tr}(\ln^2
   X)-nt+t{\rm tr}(X^{-1})-h(0)\ln\det X-h(t)\langle\langle
I_n,\exp_{X}^{-1}I_n\rangle\rangle_{X}}{({t^2+{\rm tr}(\ln^2
X))^{1/2}}}\to -\infty,$$}
   as ${\bf d_g}({\bf p}, {\bf p}_0)=({t^2+{\rm
tr}(\ln^2   X))^{1/2}}\to \infty$.
  The latter limit and the above expressions show that  hypothesis $(H_{{\bf
p}_0})$ holds true with $L_{{\bf p}_0}=-\infty.$ Now, we are in the
position to apply Theorem \ref{begle-noncompact}, i.e.,  ${\mathcal
S}_{NS}({\bf f,K})\neq \emptyset,$ while its elements $(\tilde
t,\tilde X)\in {\bf K}$ are  the solutions of the system
$$\left\{
\begin{array}{lll}
(3\tilde t^2\det \tilde X-{\rm tr}(\tilde X^{-1}))(t-\tilde t)\geq 0& {\rm \mbox{for all}} &  t\in K_1, \\
 \langle\langle -2g(\tilde t)\exp_{\tilde X}^{-1}I_n-h(\tilde t)I_n,\exp^{-1}_{\tilde
X}Y\rangle\rangle_{\tilde X}\geq 0& {\rm \mbox{for all}} & Y\in K_2.
\end{array} \right.\eqno{(S_3)}$$

Let us assume that $\tilde t=0;$ then from the first relation of
$(S_3)$ we necessarily obtain $-{\rm tr}(\tilde X^{-1})\geq 0$, a
contradiction. Thus, $\tilde t>0;$ in particular, from the first
relation of $(S_3)$ it yields that
\begin{equation}\label{ujabb-ossz-matrix}
    3\tilde t^2\det \tilde X-{\rm tr}(\tilde X^{-1})=0.
\end{equation}
Since $\tilde X\in K_2$,  the latter relation implies that
$$3\tilde t^2={\rm tr}(\tilde X^{-1})\det \tilde X^{-1}\leq {\rm tr}(\tilde X^{-1})\left(\frac{{\rm tr}(\tilde X^{-1})}{n}\right)^n\leq n.$$
This estimate gives the idea to distinguish the following two cases:

 (a) Assume $h({\sqrt\frac{n}{3}})\leq 0$. We may choose $\tilde
t={\sqrt\frac{n}{3}}$ and $\tilde X=I_n$ in $(S_3)$, taking into
account that $\langle\langle
I_n,\exp^{-1}_{I_n}Y\rangle\rangle_{I_n}=\ln \det Y\geq 0$ for every
$Y\in K_2$. Consequently, $({\sqrt\frac{n}{3}},I_n)\in \mathcal
S_{NS}({\bf f,K})$. A direct computation also shows that
$({\sqrt\frac{n}{3}},I_n)\in \mathcal S_{NE}({\bf f,K})$.

(b) Assume $h({\sqrt\frac{n}{3}})> 0$. We define the function
$j:(0,{\sqrt\frac{n}{3}}]\to [1,\infty)$ by
$j(t)=\left({\sqrt\frac{n}{3}}t^{-1}\right)^\frac{1}{n+1}.$ Since
$\lim_{t\to 0^+}j(t)=+\infty$,  $j({\sqrt\frac{n}{3}})=1$, relation
(\ref{auxiliar-egyenlte-2}), the continuity of the functions
$g,h,j$, and our assumption imply that the equation
\begin{equation}\label{ghj}
   j(t) \ln j(t)=\frac{h(t)}{2g(t)}
\end{equation}
has at least a solution $\tilde t$ with $0<\tilde
t<{\sqrt\frac{n}{3}}.$ We claim that $(\tilde t,\tilde X)=(\tilde
t,j(\tilde t) I_n)\in {\bf K}$ solves $(S_3)$, i.e., $(\tilde
t,j(\tilde t) I_n)\in {\mathcal S}_{NS}({\bf f,K})$. First of all,
we have that $\tilde X=j(\tilde t) I_n\in K_2$; indeed, ${\rm
tr}(\tilde X^{-1})=(j(\tilde t))^{-1}n< n.$ Then, a simple
calculation shows that $((\tilde t,j(\tilde t) I_n))$ verifies
(\ref{ujabb-ossz-matrix}). It remains to verify the second relation
in $(S_3);$ to complete this fact, a direct calculation and
(\ref{ghj}) give that
$$\langle\langle \exp_{\tilde X}^{-1}I_n+\frac{h(\tilde t)}{2g(\tilde t)}I_n,\exp^{-1}_{\tilde
X}Y\rangle\rangle_{\tilde X}=j(\tilde t)\left(\frac{h(\tilde
t)}{2g(\tilde t)}-j(\tilde t)\ln j(\tilde t)\right)\langle\langle
I_n,\exp^{-1}_{\tilde X}Y\rangle\rangle_{\tilde X}=0,$$ which
concludes our claim.

We also have that $(\tilde t,j(\tilde t) I_n)\in {\mathcal
S}_{NE}({\bf f,K})$.  This fact can be proved either by a direct
verification based on matrix-calculus or by observing that $X\mapsto
f_2(\tilde t,X)$ is convex (since $h(\tilde t)>0$). \hfill
$\diamondsuit$

\begin{remark}\rm
It is clear that (\ref{ghj}) can have multiple solutions which
provide distinct Nash(-Stampacchia) equilibria for ${\bf (f,K)}$.
\end{remark}


\end{example}

\begin{example}\rm
(a) Assume that $K_i$ is closed and convex in the Euclidean space
$(M_i,g_i)=(\mathbf R^{m_i},\langle\cdot,\cdot \rangle_{\mathbf
R^{m_i}})$, $i\in \{1,...,n\}$, and let ${\bf f}\in \mathcal
C_{({\bf K, U,} \mathbf R^m)}$ where $m=\sum_{i=1}^n m_i$. If
$\partial_C^\Delta {\bf f}$ is $L-$globally Lipschitz and
$\kappa$-strictly monotone on ${\bf K}\subset \mathbf R^m$, then the
function ${\bf f}$ verifies $(H_{\bf K}^{\alpha,\rho})$ with
$\alpha=\frac{\kappa}{L^2}$ and $\rho=\frac{\kappa^2}{2L^2}$. (Note
that the above facts imply that
$\kappa \leq L$, thus $0<\rho<1$.) 
Indeed, for every ${\bf p,q\in K}$ we have
that\\
\\
${\bf d}_{\bf g}^2(\exp_{\bf p}(-\alpha \partial_C^\Delta {\bf
f}({\bf p})),\exp_{\bf
 q}(-\alpha
\partial_C^\Delta {\bf f}({\bf q})))$
\begin{eqnarray*}
&=&\|{\bf p}-\alpha \partial_C^\Delta {\bf f}({\bf p})-({\bf
q}-\alpha \partial_C^\Delta {\bf f}({\bf q}))\|_{\mathbf
R^{m}}^2=\|{\bf p}-{\bf q}-(\alpha \partial_C^\Delta {\bf f}({\bf
p})-\alpha \partial_C^\Delta {\bf f}({\bf q}))\|_{\mathbf
R^{m}}^2\\&=&\|{\bf p}-{\bf q}\|_{\mathbf R^{m}}^2 -2 \alpha \langle
{\bf p}-{\bf q},\partial_C^\Delta {\bf f}({\bf p})-\partial_C^\Delta
{\bf f}({\bf q})\rangle_{\mathbf R^{m}} + \alpha^2\|
\partial_C^\Delta {\bf f}({\bf p})-\alpha \partial_C^\Delta {\bf
f}({\bf q})\|_{\mathbf R^{m}}^2\\&\leq&
(1-2\alpha\kappa+\alpha^2L^2) \|{\bf p}-{\bf q}\|_{\mathbf
R^{m}}^2=\left(1-\frac{\kappa^2}{L^2}\right) {\bf d}_{\bf g}^2({\bf
p},{\bf q})\\&\leq& (1-\rho)^2 {\bf d}_{\bf g}^2({\bf p},{\bf q}).
\end{eqnarray*}

(b) Let
$$K_1=[0,\infty),\ K_2=\{X\in M_n(\mathbf R):{\rm tr}(X)\geq 1\},$$
and the functions $f_1,f_2:K_1\times K_2 \to \mathbf R$  defined by
$$  f_1(t,X)=g(t)-ct{\rm tr}(X),\ f_2(t,X)={\rm tr}((X-h(t)A)^2).$$
Here, $g,h:K_1\to \mathbf R$ are two
 functions such that $g$ is of class $C^2$ verifying
 \begin{equation}\label{auxiliar-egyenlte-3}
0<\inf_{K_1} g''\leq \sup_{K_1} g''<\infty,
\end{equation} $h$ is $L_h-$globally
 Lipschitz, while $A\in M_n(\mathbf R)$ and $c>0$ are fixed
such that
\begin{equation}\label{utolsok-kozul}
c+L_h\sqrt{{\rm tr}(A^2)}<2\inf_{K_1} g''\  \ {\rm and}\ \
cn+2L_h\sqrt{{\rm tr}(A^2)}<4.
\end{equation}
Now, we consider the space $M_n(\mathbf R)$ endowed with the inner
product defined in (\ref{eukl-metrik-matrix}). We observe that $K_2$
is geodesic convex but  not compact in $(M_n(\mathbf
R),\langle\cdot, \cdot\rangle)$. After a natural extension of
functions $f_1(\cdot,X)$ to $\mathbf R$ and $f_2(t,\cdot)$ to the
whole $M_n(\mathbf R)$, we can state that ${\bf f}=(f_1,f_2)\in
\mathcal C_{({\bf K,U,M})},$ where ${\bf U}={\bf M}=\mathbf R\times
M_n(\mathbf R)$. 
On account of (\ref{auxiliar-egyenlte-3}), after a computation it
follows that the map
$$\partial_C^\Delta {\bf f}(t,X)=(g'(t)-c{\rm tr}(X),2(X-h(t)A))$$
is $L-$globally Lipschitz and $\kappa$-strictly monotone on ${\bf
K}$ with
$$L=\max\left(( 2\sup_{K_1} g''+8{L_h{{\rm
tr}(A^2)}})^{1/2},\left(2c^2n+8\right)^{1/2}\right)>0,$$
$$ \kappa=\min\left(\inf_{K_1} g''-\frac{c}{2}-\frac{L_h\sqrt{{\rm
tr}(A^2)}}{2},1-\frac{cn}{4}-\frac{L_h\sqrt{{\rm
tr}(A^2)}}{2}\right)>0.$$ According to (a), ${\bf f}$ verifies
$(H_{\bf K}^{\alpha,\rho})$ with $\alpha=\frac{\kappa}{L^2}$ and
$\rho=\frac{\kappa^2}{2L^2}$. On account of Theorem
\ref{dyn-system-theorem-fo}, the set of
 Nash-Stampacchia equilibrium points for
${\bf (f,K)}$ contains exactly one point $(\tilde t,\tilde X)\in
{\bf K}$ and the orbits of both
  dynamical systems $(DDS)_\alpha$ and
$(CDS)_\alpha$ exponentially converge to $(\tilde t,\tilde X)$.
 Moreover, one also has that
${\bf f}\in \mathcal K_{({\bf K,U,M})}$; thus, due to Theorem
\ref{equi-kritikus} (ii) we have that $\mathcal S_{NE}({\bf
f,K})=\mathcal S_{NS}({\bf f,K})=\{(\tilde t,\tilde X)\}.$
 \hfill $\diamondsuit$
\end{example}




\end{document}